\documentclass[12pt,a4paper]{article}
\usepackage[]{geometry}
\usepackage[mathcal]{euscript}
\usepackage{bm,url}                     
\usepackage{amsfonts}
\usepackage{amssymb}
\usepackage{amsmath}
\usepackage{todonotes}
\usepackage{amsthm}
\usepackage{setspace}
\usepackage{graphicx}               
\usepackage{natbib}                 
\usepackage{colortbl}               
\usepackage{booktabs}
\newtheorem{definition}{Definition}
\newtheorem{remark}{Remark}
\newtheorem{theorem}{Theorem}

\newtheorem{proposition}{Proposition}

\newtheorem{corollary}{Corollary}

\usepackage[english]{babel}
\usepackage{amsmath}
\usepackage{amsfonts}
\usepackage{amssymb}
\usepackage{amsthm}
\usepackage{graphicx}
\usepackage{color}
\usepackage{array}
\usepackage{mathrsfs}
\usepackage{setspace}
\usepackage{hyperref}
\definecolor{violet}{rgb}{0.7,0,0.6}

\newcommand{\eps}{\varepsilon}

\usepackage{amsfonts}
\usepackage{amssymb}
\usepackage{amsmath}
\usepackage{amsthm}

\newcommand{\X}{\mathcal{X}}

\newcommand{\ang}{\measuredangle}
\definecolor{violet}{rgb}{0.7,0,0.6}

\begin{document}
	
	\begin{center}
	\Large \bf  On estimation  of biconvex sets\\
\end{center}
\normalsize

\

\begin{center}
	Alejandro Cholaquidis$^a$, Antonio Cuevas$^b$\\
	$^a$ Universidad de la Rep\'ublica, Uruguay\\
	$^b$ Universidad Aut\'onoma de Madrid,  Espa\~na\\
\end{center}

\begin{abstract}   A set in the Euclidean plane is said to be biconvex if,  for some angle $\theta\in[0,\pi/2)$, all its sections along straight lines with   inclination angles  $\theta$ and $\theta+\pi/2$ are convex sets (i.e, empty sets or segments). Biconvexity is a natural notion with some useful applications in optimization theory. It has also be independently used, under the name of ``rectilinear convexity'', in computational geometry. We are concerned here with the problem of asymptotically reconstructing (or estimating) a biconvex set $S$ from a random sample of points drawn on $S$. By analogy with the classical convex case, one would like to define the ``biconvex hull" of the sample points as a natural estimator for $S$. However, as previously pointed out by several authors, the notion of ``hull'' for a given set $A$ (understood as the ``minimal'' set including $A$ and having the required property) has no obvious, useful translation to the biconvex case. This is in sharp contrast with the well-known elementary definition of convex hull. Thus, we have selected the most commonly accepted notion of ``biconvex hull'' (often called ``rectilinear convex hull''): we first
provide additional motivations for this definition, proving some useful relations with other convexity-related notions. Then, we prove some results 	concerning the consistent approximation  of a biconvex set $S$ and  and the corresponding biconvex hull. An analogous result is also provided for the boundaries. A method to approximate, from a sample of points on $S$, the biconvexity angle $\theta$ is also given.

\end{abstract}

\section{Introduction}\label{sec:intro}

\begin{figure}[!ht]%
	\centering
\includegraphics[scale=0.2]{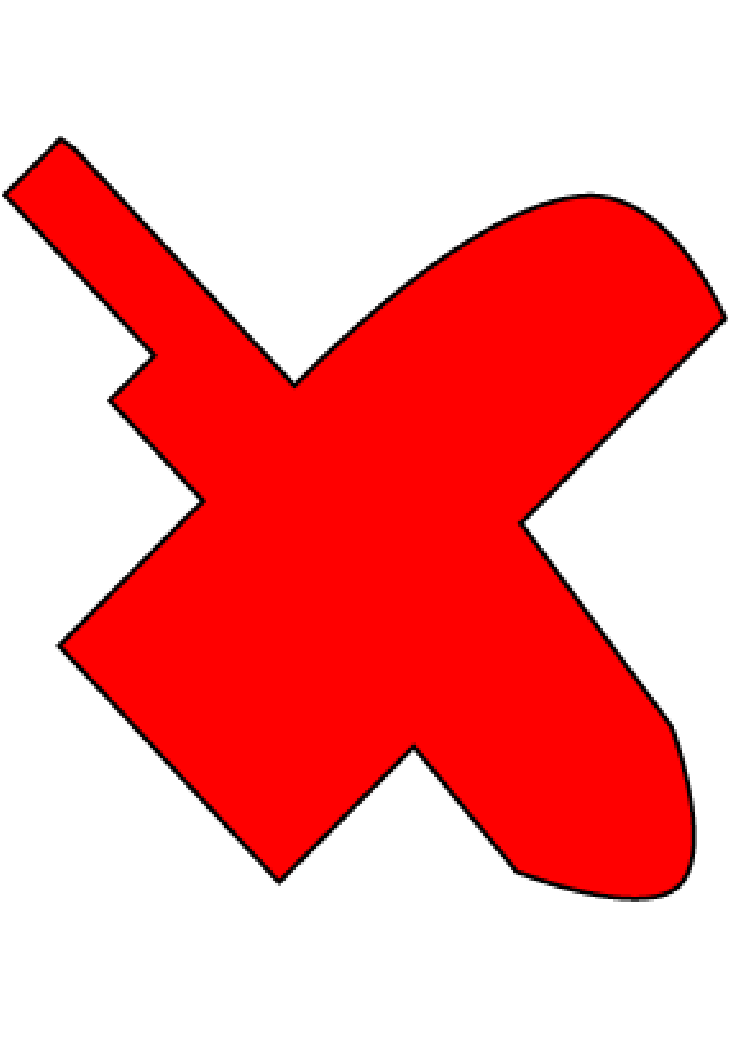} \hspace{3cm} 	\includegraphics[scale=0.5]{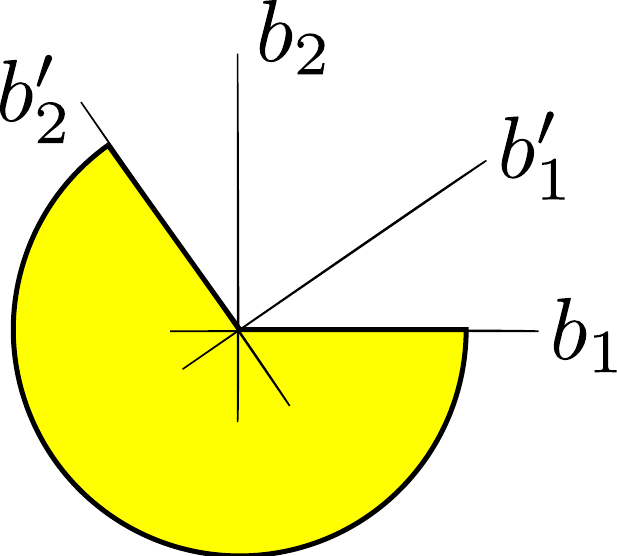} 
	\caption{On the left panel, a biconvex set. On the right panel a set that is both $(b_1,b_2)$-biconvex and $(b_1^\prime,b_2^\prime)$-biconvex.}
	\label{dibconv}
\end{figure}

 We will say that a set $S\subset \mathbb{R}^2$ is  {\bf intrinsically biconvex}, or just \textbf{biconvex}, if \textbf{there exist} two orthogonal vectors $b_1,b_2$ such that for all $\alpha,\beta\in {\mathbb R}$ the sets  $A_{\alpha}=\{\alpha b_1+ t b_2: t\in{\mathbb R}\}\cap S$ and $B_{\beta}=\{t b_1+\beta b_2: t \in {\mathbb R}\}\cap S$ are convex subsets of $\mathbb{R}^2$. When this condition is fulfilled for some given $b_1,b_2$ we will say that $S$ is biconvex with respect to (wrt) the directions $b_1$ and $b_2$.
It is clear that a set might be biconvex with respect to many different bases, see Figure \ref{dibconv}.

This notion has been used many times in the literature, for the (more restrictive) case in which the stated condition must be fulfilled for some $b_1$ and $b_2$ fixed in advance. In that case, we say  that the set $S$ is $(b_1,b_2)$-biconvex.  This concept is often expressed in terms of the inclination angle $\theta$ of the direction defined by $b_1$; we  thus can also say that the \textbf{$(b_1,b_2)$-biconvex set} $S$ is \textbf{$\theta$-biconvex}, where $\theta\in[0,\pi/2)$ is such that $b_1=(\cos(\theta),\sin(\theta))$ and $b_2=(-\sin(\theta),\cos(\theta))$. While the name $\theta$-biconvex is more convenient, we will also keep the notation $(b_1,b_2)$-biconvex for technical reasons when the reference to the biconvexity directions is useful.    The expressions \textit{(double) directional convexity}, \textit{rectilinear convexity} and \textit{restricted orientation convexity} are also sometimes used in the literature to denote this property. 
Some references are \cite{ale18} \cite{bae09}, \cite{fin98}, \cite{ott84} and \cite{raw91}.

Biconvexity is a simple extension of the classical concept of convex set. It is quite obvious that any convex set is biconvex but the converse is not true.   Such ``extended convexity'' idea (sometimes translated to functions, rather than sets) has attracted the interest of some researchers in optimization and econometrics, see \cite{aum86} and \cite{gor07} on the grounds of keeping, as much as possible, the good properties of convex functions in optimization problems.

\subsection{The set estimation point of view}\label{subsec:framework}

In the present study, we have arrived to the notion of biconvexity from a third motivation, different from computational geometry or optimization issues. Such motivation is of a statistical nature, concerning the so-called set estimation problem. The most basic version of this problem  is very simple to state: let $P_X$ be the distribution of a random variable $X$ with values in ${\mathbb R}^d$ whose support $S$ is a compact set. We aim at estimating $S$ from a random sample  $X_1,\ldots,X_n$ of independent identically distributed (iid) points drawn from $P_X$. Here the term ``estimating'' is used in the statistical sense of ``approximating as a function of the sample data''.
A consistent ``estimator'' of $S$ will be, in general, a sequence of sets,  $S_n(X_1,\ldots,X_n)$
approaching (in some suitable sense)  the set $S$ as $n$ tends to infinity.

Major applications of set estimation arise  in statistical quality control, cluster analysis, image analysis and econometrics; see the surveys by \cite{cue09} and \cite{cue10} for details. A special attention, as measured by number of citations, has deserved an application in ecology, known as \textit{home range estimation}; see, e.g.,  \cite{get04} and references therein. 

We will make no attempt to provide a complete perspective or a bibliography on set estimation. The previous remarks aim only at establishing the setting in which the present study must be included, thus providing some insight to interpret our results. In addition to be above mentioned survey papers, we refer also to \cite{chola:14}, \cite{aar16}, \cite{che17} and references therein, for more recent contributions on this and other closely related topics.

\subsection{The plan and contributions of this work} \label{subsec:plan}

We aim at exploring the applicability of the notion of biconvexity  in the above mentioned statistical problem of reconstructing, from a random sample of points, a two-dimensional compact set $S$.

In Section \ref{sec:definitions} we will introduce some notation and auxiliary definitions. 

In Section \ref{sec:properties} we will relate the concept  of biconvex set with the notion of ``lighthouse set'' previously analyzed in \cite{chola:14}.

In Section \ref{sec:hull} we will consider the problem of estimating an unknown biconvex set in ${\mathbb R}^2$ from a random sample of points ${\mathcal X}_n=\{X_1,\ldots,X_n\}$ whose distribution has support $S$.  We propose to estimate $S$ using a biconvex hull, ${\mathbb B}({\mathcal X}_n)$, which (from a completely different point of view) has been previously considered in the literature on computational geometry; see e.g., \cite{bae09}. In particular, we will prove (in Theorem \ref{thhull}) the statistical consistency, as well as convergence rates, for the estimator ${\mathbb B}({\mathcal X}_n)$ with respect to the Hausdorff metric and the ``distance in measure'' commonly used in set estimation problems. An additional result concerning the estimation of the true biconvexity angle will be also proved in Theorem \ref{th:anglechoice}. 

Some numerical illustrations are included in Section \ref{sec:sim}.

 Overall, the main achievement of this paper is to analyze, from the statistical point of view, the class of biconvex sets in the plane. We show that, under quite reasonable additional regularity properties, these sets can be estimated with a reasonable simplicity.  Also, from the point of view of computational geometry, we provide some additional compelling reasons (see Theorem \ref{thhull} below) for the use of the ``rectilinear convex hull'' (see \cite{bae09} and references therein) as a natural notion of ``biconvex hull'' of a finite sample of points. 

\section{Some notation and definitions}\label{sec:definitions}

We consider $\mathbb{R}^2$ endowed with the Euclidean norm $\|\cdot\|$. The closed ball of radius $r$ centred at $x$ is denoted by $B(x,r)$. The interior of the ball is denoted by $\mathring{B}(x,r)$. With a slight abuse of notation, if $S\subset \mathbb{R}^2$, we will denote the $r$-parallel set by $B(S,r)=\cup_{s\in S} B(s,r)$.  The two-dimensional Lebesgue measure will be denoted $\mu$ and $\omega_2=\mu(B(0,1))$. For $\eps>0$ and $A\subset \mathbb{R}^2$, we define $A\ominus B(0,\eps)=\{x:B(x,\epsilon)\subset A\}$. The distance from a point $x$ to $S$ is denoted by $d(x,S)$, i.e: $d(x,S)=\inf\{\|x-s\|:s\in S\}$.  If $S\subset\mathbb{R}^2$, $\partial S$, $int(S)$ (or $\mathring{S}$), $S^c$, $\overline{S}$ stand for the boundary, interior, complement, and topological closure of $S$, respectively. Given two points $p_1$ and $p_2$ we denote $\overline{p_1p_2}$ the closed segment joining $p_1$ and $p_2$. Given a coordinate system $[b_1,b_2]$ we denote $R$ the counter clockwise rotation of angle $\pi/2$ with center at $(0,0)$ and $R_\theta$ the clockwise rotation of angle $\theta\in [0,\pi)$. If $\theta\in (-\pi,0]$, $R_\theta$ is the $\theta$-counter clockwise rotation. Given two vectors $v_1,v_2$, we define  $\ang(v_1,v_2)=\arccos(\langle v_1,v_2 \rangle /(\|v_1\|\|v_2\|))\in [0,\pi]$.
For $\xi\in{\mathbb R}^2$,   $R^i(\xi)$   will represent  (for $i=1,2,3$)  the counter-clockwise $\pi/2$-rotation of $R^{i-1}(\xi)$, with $R^0(\xi)=\xi$.  The canonical basis in $\mathbb{R}^2$ will be denoted $\{e_1,e_2\}$.

\

\noindent \textit{Lighthouses}

An infinite (open) cone with vertex $x$ is defined by
$$C_{\rho,\xi}(x):=  \left\{z\in {\mathbb R}^2,\, z\neq x: \Big \langle \xi,\,\frac{z-x}{\|z-x\|} \Big \rangle >\cos(\rho/2)\right\},$$
where $\xi\in{\mathbb R}^2$, with $\Vert \xi\Vert=1$ is the direction of the cone axis and $\rho\in(0,\pi]$ is the ``opening angle''.  So, in particular, $C_{R(\xi)}(x)$ will stand in what follows for the cone with vertex $x$, $\rho=\pi/2$ and axis $R(\xi)$.

The following notions of ``lighthouse-sets''  were introduced in \cite{chola:14}. We will use them here as auxiliary tools, conceptually related to the concept of biconvex set we established at the very beginning of the paper. 

\begin{definition}\label{lightcomp}
	Given a set  $S\subset\mathbb{R}^2$ and an opening angle $\rho\in(0,\pi]$, we define the $\rho$-\textbf {lighthouse hull by complements} of $S$ by
	\begin{equation}\label{eq:ccc}
	{\mathbb C}_\rho(S)=\bigcap_{\{y,\xi\in{\mathbb R}^2:\Vert \xi\Vert=1:\,C_{\rho,\xi}(y)\cap S=\emptyset\}}\left(C_{\rho,\xi}(y)\right)^c.
	\end{equation}
	A  set $S\subset\mathbb{R}^2$ is said to be a \textbf{$\rho$-lighthouse by complements} when $S={\mathbb C}_\rho(S)$.
	
	Finally, $S$ is said to be a \textbf{$\rho$-lighthouse set} if for each $x\in\partial S$ there exists an open cone $C_{\rho,\xi}(x)$ with vertex at $x$ such that $C_{\rho,\xi}(x)\cap S=\emptyset$.  
\end{definition}

\begin{figure}[!ht]%
	\centering
	\includegraphics[scale=0.5]{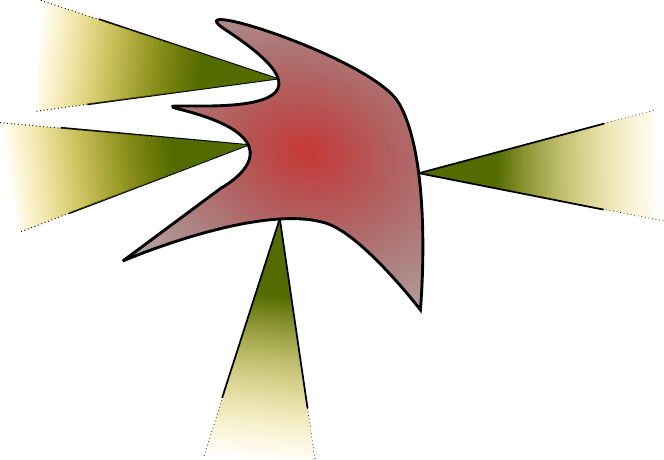} 

	\caption{The supporting cone property of lighthouse sets}
	\label{dib:lighthouse}
\end{figure}

In other words, a set  $S$ is a $\rho$-lighthouse by complements if and only if $S$ can be expressed as the intersection of the complements of all open cones of type $C_{\rho,\xi}$ that are disjoint with $S$. 
Also, $S$ is just a $\rho$-lighthouse set if for every boundary point $x$ there is a ``supporting cone'' with vertex at $x$  that is completely included in $S^c$. Such supporting cone could be seen, in intuitive terms, as an undisturbed ``beam of light'' projected outside $S$ from any point of $\partial S$.  See Figure \ref{dib:lighthouse}.

In \cite{cholatesis} Proposition 3.6 d), it is  proved that if $S$ is a $\rho$-lighthouse by complements, then it is a $\rho$-lighthouse set as well. The converse implication is not true in general; see Figure 3.2 in \cite{cholatesis}.

It is clear that a compact set $S$ is convex if and only if $S$ is 
a $\pi$-lighthouse by complements, according to Definition \ref{lightcomp}. 

 Also, if the ``volume elements" $C_\rho(x)$ are replaced with open balls $B(x,\alpha)$ in Definition \ref{lightcomp} we get a related concept, often called \textit{$\alpha$-convexity}. Finally the above mentioned ``cone supporting property'' boils down to the so-called \textit{outer $\alpha$-rolling ball property} when $C_{\rho,\xi}(x)$ is replaced with $B(x,\alpha)$; see \cite{cue12}, \cite{ari18} for additional information on $\alpha$-convexity and rolling properties. 

\

\noindent \textit{Metrics between sets, boundary measure}

The performance of a set estimator 
is usually evaluated through the Hausdorff distance (\ref{dh}) and the distance in measure (\ref{dmu}) given below. 
The distance in measure takes the mass of the symmetric difference into account while the Hausdorff distance measures, in some sense, the difference of the shapes.

Let $A,C\subset \mathbb{R}^2$ be non-empty compact sets. The Hausdorff (or Hausdorff-Pompeiu) distance between $A$ and $C$ is defined as 
\begin{equation}\label{dh}
d_H(A,C)=\max\Big\{\max_{a\in A}d(a,C), \ \max_{c\in C}d(c,A)\Big\}.
\end{equation}

If $\nu$ is a Borel measure, the distance in measure between  $A$ and $C$ is defined as 
\begin{equation}\label{dmu}
d_\nu(A,C)=\nu(A\triangle C),
\end{equation}
where $\triangle$ denotes symmetric difference.

The following notion of ``boundary measure'' is quite popular in geometric measure theory as a simpler alternative to the more sophisticated notion of Hausdorff measure. See \cite{amb08}, and references therein, for the geometric aspects of this concept. See \cite{cue12} and \cite{cue18} for some statistical applications.

\begin{definition} \label{minkcont} Let $S\subset\mathbb{R}^2$ be a compact set.  The Minkowski content of $\partial S\subset \mathbb{R}^2$  is given by,
	$$L_0(\partial S)=\lim_{\epsilon\rightarrow 0} \frac{\mu\big(B(\partial S,\epsilon)\big)}{2\epsilon},$$
	provided that the limit exists and it is finite.
\end{definition}

Note that the existence and finiteness of $L_0(\partial S)$ is itself a regularity condition on the set $S$. We will need such condition later to establish some of our results (in particular Theorems \ref{thhull} and \ref{th:anglechoice} and Corollary \ref{cr_revisado}).

\section{The ``lighthouse properties'' of regular biconvex sets} \label{sec:properties}

The purpose of this section is to show that biconvexity is essentially equivalent (under some regularity conditions,  which will be shown to be necessary in order to get the equivalence) to a restricted version of the lighthouse properties, established in Definition \ref{lightcomp}, in which $\rho=\pi/2$ and the direction of the cone axes are fixed   up to a $\pi/2$-rotation.   The formal statement is as follows. 

\

\begin{theorem} \label{p1} Let $S\subset \mathbb{R}^2$ be a closed set such that $\partial S$ is path-connected, and $S=\overline{int(S)}$. Then, $S$ is biconvex with respect to the orthonormal vectors $b_1,b_2$ if and only if for all $x\in \partial S$, there exists $i\in\{0,1,2,3\}$, such that $C_{R^i(\xi)}(x)\cap S=\emptyset$, where $\xi=(b_1+b_2)/\|b_1+b_2\|$, $R^0(\xi)=\xi$ and, for $i>0$,  $R^i(\xi)$  stands by the counter-clockwise $\pi/2$-rotation of $R^{i-1}(\xi)$.  
	\begin{proof} Assume for simplicity that $\{b_1,b_2\}$ is the canonical basis $\{e_1,e_2\}$, so that $b_2$ corresponds to the vertical direction. Let us assume that $S$ is biconvex wrt $e_1,e_2$, we will prove that is a $\pi/2$-lighthouse with $\xi=(e_1+e_2)/\|e_1+e_2\|$. If this is not the case there exists $x\in \partial S$ and $z_i=(\alpha_i,\beta_i)\in C_{R^i(\xi)}(x)\cap S$ for $i=0,1,2,3$, where $(\alpha_i,\beta_i)$ are the coordinates with respect to $\{e_1,e_2\}$, we may assume that $x=(0,0)$ (otherwise we could consider a translation of $S$).  Since $\partial S$ is path connected, $S$ is path-connected.   Let $\gamma_i$ be a path   in $S$  connecting $z_{i}$ with $(0,0)$ for $i=0,1,2,3$.
	Let $t=(t_1,t_2)\in S^c$ with $t_1<(1/2)\min_i \alpha_i$ and $t_2<(1/2)\min_i \beta_i$. Suppose none of the four paths meets $\{(t_1,t):t>t_2\}\cap \{(t,t_2):t>t_1\}$ then, since $t\in S^c$, non of the paths meets $\{(t_1,t):t\geq t_2\}\cap \{(t,t_2):t\geq t_1\}$, which is not possible for $\gamma_0$. Reasoning in the same way with the other four cones centred at $t$ we get a contradiction.

		To prove the other implication let us assume that $S$ is $\pi/2$-lighthouse with  possible  axes $\xi,R(\xi),$ $R^2(\xi),$ $R^3(\xi)$ were $\xi=(e_1+e_2)/\|e_1+e_2\|$,  but $S$ is not biconvex wrt $e_1,e_2$. We have two possibilities:
		\begin{itemize}
			\item[1)]  there exist $z_1=\alpha_1e_1+\beta_1e_2\in S$ and $z_2=\alpha_2e_1+\beta_1e_2\in S$ such that the  ``horizontal'' set  $\overline{z_1z_2}\cap S$ is not convex,
			\item[2)] there exist $z_1=\alpha_1e_1+\beta_1e_2\in S$ and $z_2=\alpha_1e_1+\beta_2e_2\in S$ such that the  ``vertical'' set  $\overline{z_1z_2}\cap S$ is not convex. 
		\end{itemize}
	We will consider the first case,  as the  second one is analogous. From the non-convexity of $\overline{z_1z_2}\cap S$ there exist  $p_1\neq p_2$ such that $p_1\in \overline{z_1z_2}\cap \partial S$, $p_2\in \overline{z_1z_2}\cap \partial S$ and $\overline{p_1p_2}\cap S=\{p_1,p_2\}$.   To simplify the notation we will assume (without loss of generality) that $p_1=(0,0)$ and $p_2=(a,0)$ for some $a>0$. 
	Let us consider any curve $\gamma$,  included in $\partial S$,   joining $p_1$ with $p_2$. We may assume  without loss of generality that $\gamma$ does not have auto-intersections. Denote by $\rho(t)=(\rho_1(t),\rho_2(t))$ for $t\in[0,1]$ the part of the curve $\gamma$ whose abscissa is always between $0$ and $a$, see Figure \ref{dib2}, that is $0\leq \rho_1(t)\leq a$. Since $\rho$ does not intersect the open segment $\overline{p_1,p_2}\setminus \{p_1,p_2\}$, we have two possibilities $\rho_2(t)>0$ for all $t$ or $\rho_2(t)<0$ for all $t$ (recall that $p_1=(0,0)$ and $p_2=(a,0)$), assume that we are in this last case (the first one is analogous).
	Let us denote $H_\rho=\{(x,y): 0\leq x\leq a,  \rho_2(t)\leq y\leq 0\ \ \forall t\in [0,1]\}$, see Figure \ref{dib2}.
  Let us prove that $\mbox{\rm int}(S)\cap H_\rho=\emptyset$. Suppose by contradiction that there exists $l=(l_1,l_2)\in\mbox{\rm int}(S)\cap H_\rho$, let $\delta>0$ such that $B:=B(l,\delta)\subset S$,  then there exists $v=(v_1,v_2)\in \partial S\cap H_\rho$ in the perpendicular line to $\overline{p_1p_2}$ passing through $l$ and with $l_2<v_2<0$. But this imply that the lighthouse property (with  possible  axes $\xi,R(\xi),$ $R^2(\xi),$ $R^3(\xi)$) fails to be fulfilled in $v$ since $p_1\in C_{R(\xi)}(v)$,  $p_2\in C_{\xi}(v)$, $C_{R^2(\xi)}(v)\cap B\neq \emptyset$ and $C_{R^3(\xi)}(v)\cap B\neq \emptyset$. 
  A completely similar reasoning leads to $\text{int}(S)\cap G_\rho=\emptyset$ and $G_\rho=\{(x,y):0\leq x\leq a, y< \rho_2(t)\,\forall t\in [0,1]\}.$

  Finally we have proved, $\text{int}(S)\cap H_\rho=\emptyset$ and $\text{int}(S)\cap G_\rho=\emptyset$ which contradicts $S=\overline{\text{int}(S)}$.
\begin{figure}[!ht]%
	\centering
	\includegraphics[scale=0.3]{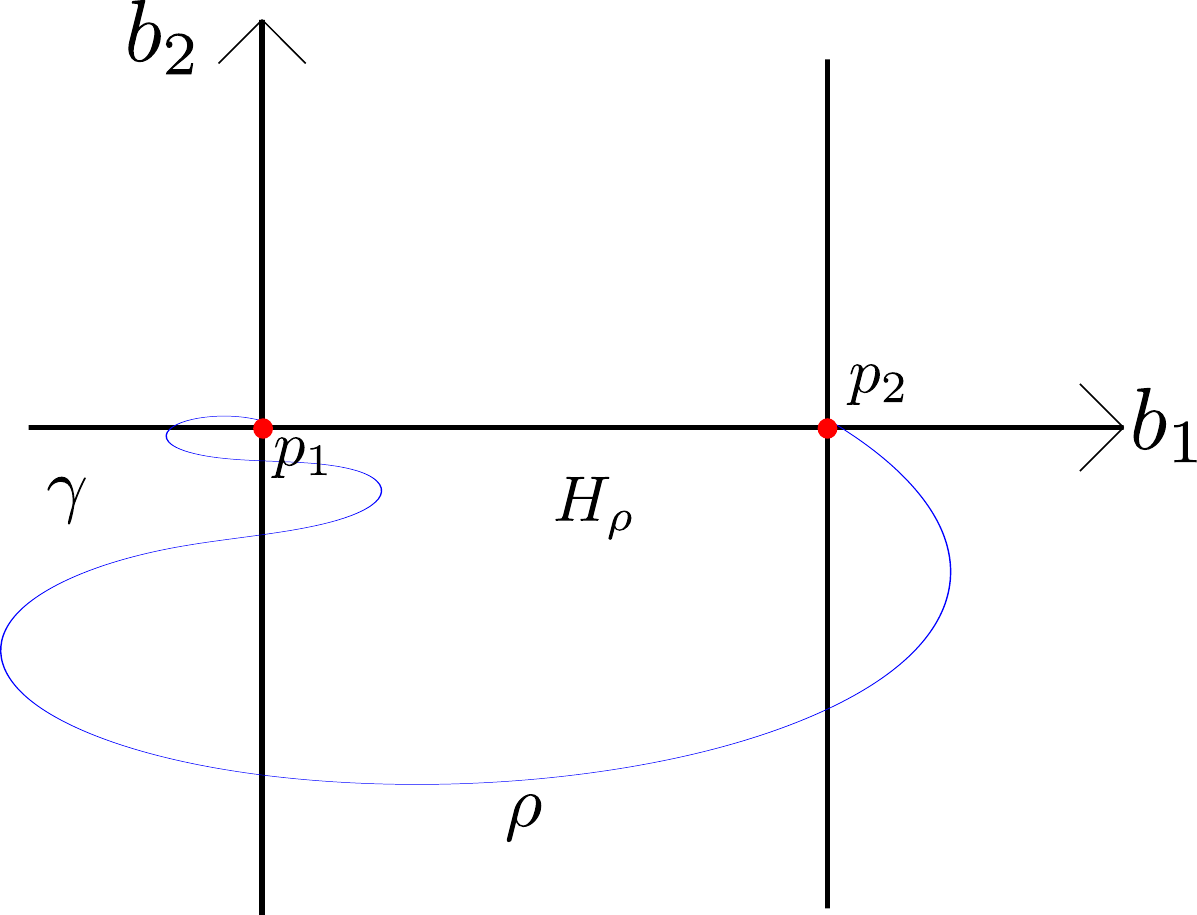} \ \ 	
	\caption{The curve $\gamma$ joins $p_1=(0,0)$ and $p_2=(a,0)$ and the curve $\rho=(\rho_1,\rho_2)$ is the part of $\gamma$ such that $0\leq \rho_1\leq a$.}
	\label{dib2}
	
\end{figure}

		\end{proof}
\end{theorem}
 
\begin{remark} Both hypothesis ($\partial S$ is path-connected, and $S=\overline{int(S)}$) are necessary in order to get the equivalence. For example $\partial B(0,1)$   fulfils the lighthouse property   but it is clearly not biconvex (observe that it is path-connected). Also $\overline{B((0,0),1/4)\cup B((1,1),1/4)}$ is biconvex but it is not  a lighthouse set   when we restrict the axes of the cones to be given by some $\xi$ and $R^i(\xi)$ for $i=1,2,3$.
\end{remark}
 
The following result provides some additional insights on the geometric nature of compact biconvex sets. We prove that, under mild regularity assumptions, these sets can be expressed as the intersection of the complements of a family of open quadrants. This result will be useful later in order to study the estimation of a biconvex set from a random sample.
 
\begin{theorem} \label{th:p2} Let $S\subset \mathbb{R}^2$ be a compact set such that $\partial S$ is path-connected, and $S=\overline{\mbox{\rm int}(S)}$.   Then,   $S$ is biconvex wrt $b_1,b_2$ if and only if   there exists $\xi$ with $\|\xi\|=1$, such that for all $y \in S^c$, there exists $i\in \{0,1,2,3\}$ and $z$ such that $y \in C_{\pi/2,R^i(\xi)}(z)\subset S^c$.
	\begin{proof} If we assume that the set is $\pi/2$-lighthouse by complements then the biconvexity follows  from  Theorem \ref{p1} together with the fact that the ``lighthouse by complements'' property implies the plain lighthouse condition (see Definition \ref{lightcomp}).  Now, to prove the other implication let us assume by contradiction that $S$ is biconvex wrt $b_1,b_2$, but not a $\pi/2$-lighthouse by complements (where the axes are $\xi,R(\xi),R^2(\xi),R^3(\xi)$ and $\xi=(b_1+b_2)/\|b_1+b_2\|$).    Then there must be a point  $t$ which ``cannot be separated from $S$'' using quadrants with the prescribed axes. In more precise terms,  there   exist   $t\in S^c$ and $z_0,\dots,z_3\in S$ such that  $z_i\in \overline{C_{R^i(\xi)}(t)}\cap S$   for $i=0,1,2,3$ and $\xi=(b_1+b_2)/\|b_1+b_2\|$. Let $\gamma_i\subset S$ joining $z_i$ with $z_{i+1}$ where $z_4=z_0$. 
		Let $\delta>0$ such that $B(t,\delta)\subset S^c$. In what follows the coordinates are in the axes $b_1,b_2$ and $t=(0,0)$. If there exists $0<y_1$ and  $y_2<0$ such that $\{(0,y_1),\, (0,y_2)\}\subset\gamma_1$, from the biconvexity it follows that $t\in S$, see Figure \ref{dibcomp} left. Clearly the same holds for $\gamma_3$. Reasoning in the same way, if there exists $0<x_1$ and $x_2<0$ such that  $\{(0,x_1),\,(0,x_2)\}\subset\gamma_0$  then $t\in S$, and the same holds for $\gamma_2$.  Finally the only other possible configuration is shown in Figure \ref{dibcomp} right, which also leads to $t\in S$ since $t$ is in the middle of a vertical (or horizontal) segment with extremes  in $S$.
	\begin{figure}[!ht]%
	\centering
	\includegraphics[scale=0.3]{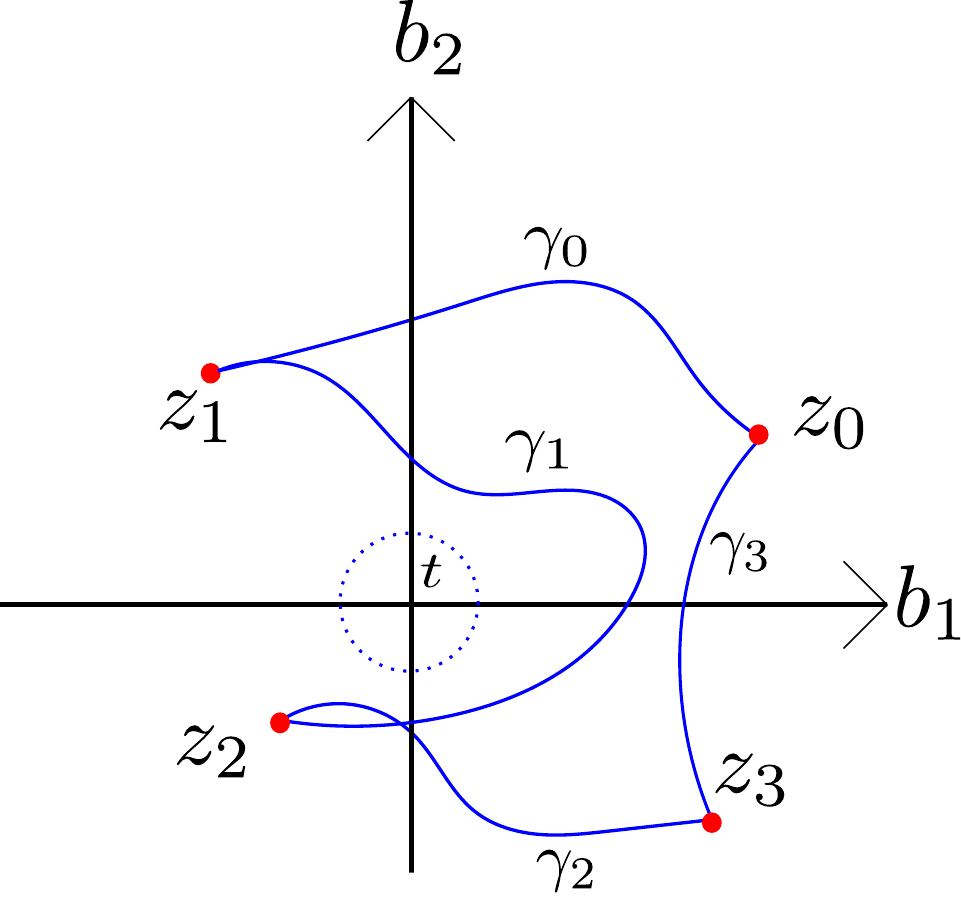} \ \ 	
	\includegraphics[scale=0.3]{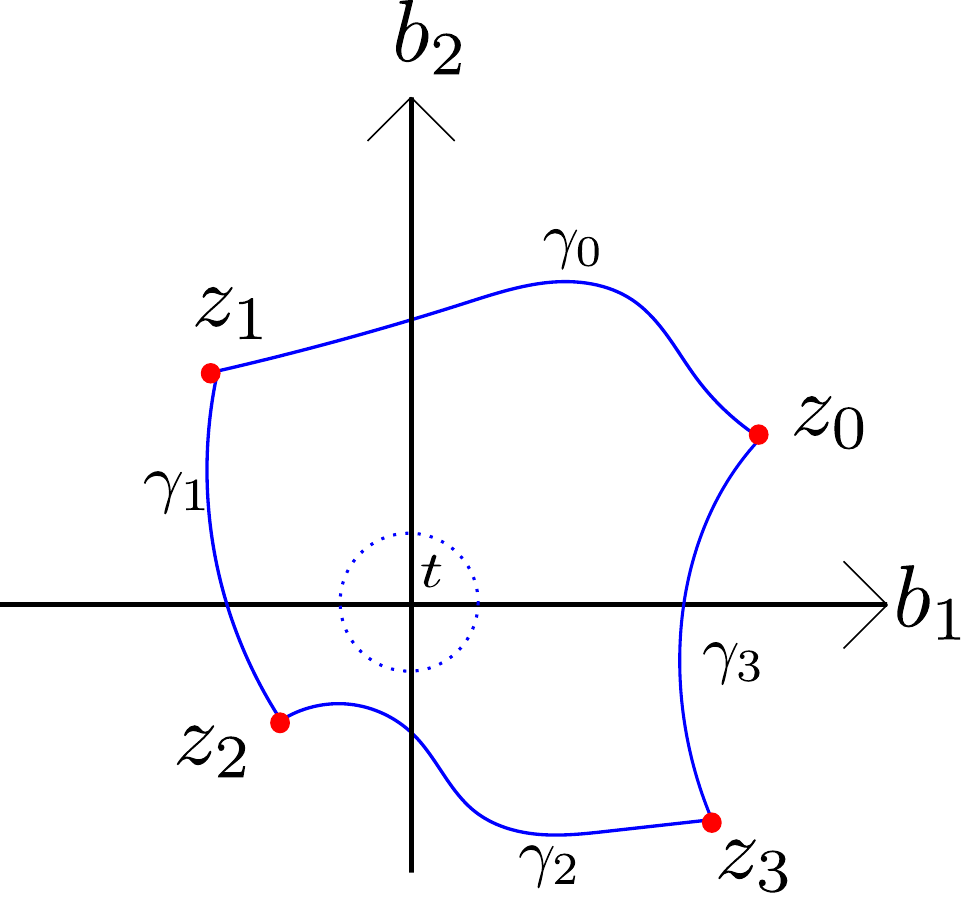}	
	\caption{Left: If there exists $0<y_1$ and $y_2<0$ such that $\{(0,y_1),\ (0,y_2)\}\subset\gamma_1$  then $t\in S$. Right:  no path $\gamma_i$,  (for $i=1,3$) meets the vertical axis both above and below $t\in S$; also, no path $\gamma_i$,  (for $i=0,2$) meets the horizontal axis both on the left and on the right of $t$.}
	\label{dibcomp}
\end{figure}
\end{proof}
\end{theorem}


\begin{remark} Theorems \ref{p1} and \ref{th:p2} imply that the class of compact \textbf{$\pi/2$-lighthouses by complements}  agrees with the class of compact   \textbf{$\pi/2$-lighthouses}, when we restrict ourselves to sets whose boundary is path connected and fulfil $S=\overline{\mbox{\rm int}(S)}$. This coincidence does not hold in general, even if we restrict ourselves to compact sets fulfilling $S=\overline{int(S)}$; see \cite{chola:14}.  In Theorem \ref{th:p2} the axe $\xi$ is not necessarily unique, and the index $i$ in general depends on $y$.
\end{remark}

\subsection{On the angles of biconvexity}
 Instead of considering the angles  as points in $[0,\pi/2]$ we will take the quotient space $[0,\pi/2]/\sim$ with the quotient topology where we identify $0=\pi/2$.  This is equivalent to view the angles as points in ${\mathbb S}^{\pi/2}:=\partial B((0,0),1/4)$; the choice of the  radius $1/4$ allows for an additional simple identification of every point $P$ in ${\mathbb S}^{\pi/2}$ with the length of the counter-clockwise arc from the point $(1/4,0)$ to $P$. 
For convenience, we will use such identification in what follows. The next proposition states that the set of  biconvexity angles is an arc (which could reduce to just a single point) in ${\mathbb S}^{\pi/2}$.

\begin{proposition} \label{propconv} Let $S\subset \mathbb{R}^2$ be a compact, biconvex set such that $\partial S$ is path-connected, and $S=\overline{int(S)}$. Then, the set of angles is an arc in ${\mathbb S}^{\pi/2}$.

\begin{proof} Since $S$ is biconvex there must be at least an angle $\theta_0$ for which the condition of biconvexity is fulfilled. If there is no other value of $\theta$ for which $S$ is $\theta$-biconvex, then the proof is concluded (in that case the arc would reduce to $\{\theta_0\}$). Otherwise, we can take $0\leq \theta_0<\theta_1<\pi/2$ such that $S$ is biconvex in the directions determined by two angles $\theta_0$ and $\theta_1$. We can assume that $\theta_0=0$, otherwise take as canonical directions the ones determined by $\theta_0$.  For $j=0,1$, let us denote by  $\xi_j=(\cos(\theta_j+\pi/4),\sin(\theta_j+\pi/4))$. 
 From Theorem \ref{p1}, $S$ is a $\pi/2$-lighthouse by considering only cones with axes given by $\xi_0$, $\xi_1$ with $\|\xi_0\|=1$, $\|\xi_1\|=1$ or any rotation of angle $\pi/2$ of $\xi_0$ or $\xi_1$. Since $\theta_0\neq \theta_1$, $\xi_1\neq R^i(\xi_0)$. 
	
Let us consider $x\in \partial S$. We will prove that we have only two possibilities (we will refer to them as Case 1 and Case 2),
	\begin{itemize}
		\item[1)] there  exist  $i,j$ depending on $x$ such that $C_{R^j(\xi_0)}(x)\cap S=\emptyset$, $C_{R^i(\xi_1)}(x)\cap S=\emptyset$ and $C_{R^i(\xi_0)}(x)\cap C_{R^j(\xi_1)}(x)\neq \emptyset$. In this case $C_{\psi}(x)\cap S=\emptyset$ for any $\psi$ in the cone determined by $R^i(\xi_0)$ and $R^j(\xi_1)$. Observe that if all the points $x\in \partial S$ are in this   Case 1  then   by Theorem 1 $S$ is lighthouse where the possible axes are $\psi$ in the aforementioned cone, (and the four $\pi/2$ rotations of $\psi$). Since we assume $\theta_0=0$, this  implies   (again by Theorem 1) that the set $S$ is biconvex wrt $\theta$,   for all $\theta\in [0,\theta_1]$.  So, the set of convexity angles is connected and therefore an interval (i.e. an arc, when viewed in ${\mathbb S}^{\pi/2}$). 
		
	    \item[2)] for some $i\in\{0,1,2,3\}$ and $j\in\{0,1\}$, $\big(C_{R^i(\xi_j)}(x)\cup C_{R^{i+1}(\xi_j)}(x)
	    \big)\cap S=\emptyset$; here $C_{R^4(\xi_j)}(x)$ stands for $C_{R^0(\xi_j)}(x)$. In this case there is a half-space $H$ not meeting $S$ such that $x\in \partial H$  (such half-space would be the topological closure of $C_{R^i(\xi_j)}(x)\cup C_{R^{i+1}(\xi_j)}(x)$); recall that we are considering open cones but, still, we cannot have points of $S$, apart form $x$, in the common half-line boundary between $C_{R^i(\xi_j)}(x)$  and $C_{R^{i+1}(\xi_j)}(x))$ due to the assumption $S=\overline{int(S)}$. If all the points $x\in \partial S$ are in this   Case 2  the set $S$ is convex and therefore biconvex for all $\theta\in [0,\pi/2]$.
 \end{itemize}

\

Suppose that  for $x\in\partial S$  we are not in Case 2, then for all $i=0,1,2,3$ and $j=0,1$, 
 \begin{equation} \label{case2}
\big(C_{R^i(\xi_j)}(x)\cup C_{R^{i+1}(\xi_j)}(x)\big)\cap S\neq\emptyset,
\end{equation} 
 we will prove that this implies that we are in Case 1.
Since $S$ is $\pi/2$-lighthouse with axes $R^i(\xi_0)$ and $R^k(\xi_1)$ for some $i,k=0,\dots,3$, there exist  $i$ and $k$ such that  $C_{R^i(\xi_0)}(x)\cap S= \emptyset$ and $C_{R^k(\xi_2)}(x)\cap S= \emptyset$.   Let us take, for example, the case $i=0$ and $k=1$ as in Figure \ref{figconv}. 
We will prove that $C_{R^2(\xi_0)}(x)\cap S=\emptyset$ and then we are in Case 1 because $C_{R^2(\xi_0)}(x)\cap C_{R(\xi_1)}(x)\neq \emptyset.$

Assume by contradiction this is not true.  Let $E_0$ be the coordinate system determined by $\theta_0$ with $x=(0,0)$, where $\xi_0$ has positive coordinates; and $E_1$ determined by $(0,0)$ and $\theta_1$, where $\xi_1$ has positive coordinates.  Then  $C_{R^2(\xi_0)}(x)\cap S\neq \emptyset$  so that  there must exist some $(b_1,b_2)\in C_{R^{2}(\xi_0)}(x)\cap S$, where  $(b_1,b_2)$ are the coordinates in $E_0$, since $C_{R(\xi_1)}(x)\cap S=\emptyset$ we have that $b_1,b_2<0$. Observe that the coordinates of this point wrt $E_1$ must be both negative also. 
Let us denote $(a_1,a_2)\in S\cap C_{R(\xi_0)}(x)$ where $(a_1,a_2)$ are the coordinates in $E_1$ (see Figure \ref{figconv}); note that such $(a_1,a_2)$ must exist since we are assuming that we are not in Case 2. Let us assume first that $a_1\leq b_1$ as in Figure \ref{figconv}. Let $\gamma\subset S$ be a curve joining $(a_1,a_2)$ with $x$, in this case there exists  $p\in \{(b_1,t):\,t\in \mathbb{R}\}\cap \gamma$. Since $\gamma\subset S$ and $S$ is biconvex $\overline{(b_1,b_2),p}\subset S$. Since $a_1<0$ then $\overline{(b_1,b_2),p}\cap C_{R(\xi_1)}(x)\neq \emptyset$, which  contradicts  that $C_{R(\xi_1)}(x)\cap S=\emptyset$. If $a_1> b_1$ a similar contradiction is obtained by considering a curve $\gamma\subset S$ joining $(b_1,b_2)$ with $x$ and a line $r=\{(a_1,t):t\in \mathbb{R}\}$.
We know that   $C_{R^2(\xi_1)}(x)\cap S= \emptyset$ then $C_{R^2(\xi_0)}(x)\cap C_{R(\xi_1)}(x)\cap S=\emptyset$.

We have thus obtained that $[(C_{R^2(\xi_0)}(x)\cap C_{R^2(\xi_1)}(x))\cup C_{R(\xi_1)}(x)]\cap S=\emptyset$.
Now, note that the set $C_{R^2(\xi_0)}(x)$ can be expressed as a union of three sets: the first one is $[(C_{R^2(\xi_0)}(x)\cap C_{R^2(\xi_1)}(x)]$; the second one is $C_{R^2(\xi_0)}(x)\cap C_{R(\xi_1)}(x)$: none of these sets intersects $S$, as we have proved  $[(C_{R^2(\xi_0)}(x)\cap C_{R^2(\xi_1)}(x))\cup C_{R(\xi_1)}(x)]\cap S=\emptyset$. The third set is the half-line $\{(0,t):t<0\}\cap S$ where $(0,t)$ are the coordinates in $E_1$; but again we have $\{(0,t):t<0\}\cap S=\emptyset$, as a consequence of the assumption $S=\overline{int(S)}$. It follows that $C_{R^2(\xi_0)}(x)\cap S=\emptyset$ so that we are in Case 1 (recall that $C_{R(\xi_1)}(x)\cap S=\emptyset$).  This conclude the proof that we are either in case 1 or case 2.

To conclude the proof of the Theorem recall that if  all the points $x\in \partial S$ are in Case 2  then $S$ is convex, and then is biconvex for all $\theta\in [0,\pi/2]$. If all points are as in case 1 then $S$ is biconvex for all $\theta \in [0,\theta_1]$. Finally if there are points in Case 1 and points in Case 2,   we would also have that $S$ is biconvex wrt $\theta$ for all $\theta\in [0,\theta_2]$, as the points in Case 2 do not introduce any restriction on the lighthouse axes.

\begin{figure}[!ht]
		\centering
	\includegraphics[scale=0.7]{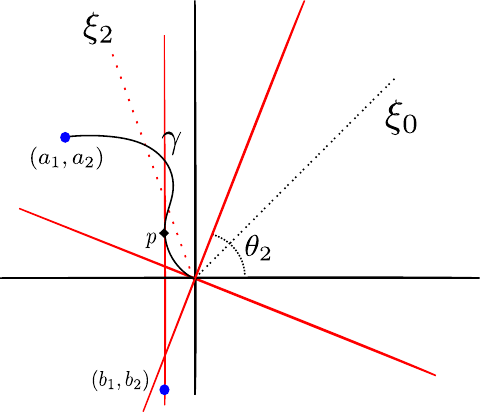} \ \ \ \ 
	\caption{The case  $a_1\leq b_1$.}
		\label{figconv}
	\end{figure}
\end{proof}
\end{proposition}

\begin{figure}[!ht]
		\centering
		\includegraphics[scale=0.5]{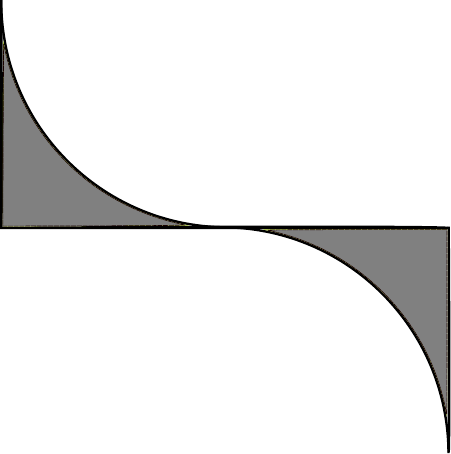} \ \ \ \ 
		\caption{$\theta_0=0$ is the only biconvexity direction}
		\label{hypograph}
	\end{figure}

\begin{remark}  Proposition  \ref{propconv} does not prove that the set of angles is always a proper non-degenerate arc. This is not true in general, as it can be seen in Figure \ref{hypograph}. The set shown is the union of two sets which are obtained by rotation and translation of the hypograph of the function $f(x)=\sqrt{1-x^2}$ for $0\leq x\leq 1$.
This set fulfils all the conditions of Proposition  \ref{propconv}; however, it is clear that $\theta_0=0$ is the only biconvexity direction in this case.
\end{remark}

\section{Statistical estimation of $\theta$-biconvex sets} \label{sec:hull}

We now consider the statistical problem of estimating a biconvex, path connected, compact set $S\subset{\mathbb R}^2$ from a sample $\mathcal{X}_n=\{X_1,\dots,X_n\}$ drawn from  a  distribution $P_X$ whose support is $S$. 

By analogy with other similar problems, based on convexity type assumptions on $S$ (such as convexity or $\alpha$-convexity; see \cite{cue12}) one would be tempted to estimate $S$ using the ``biconvex hull of ${\mathcal X}_n$'' that is, the intersection of all biconvex sets containing ${\mathcal X}_n$. However the biconvex hull of  ${\mathcal X}_n$ will be, in most cases, the sample ${\mathcal X}_n$ itself, since typically ${\mathcal X}_n$ will be biconvex   with respect to some  orthonormal directions $b_1$ and $b_2$ given in advance.   For example, this will happen with probability one  whenever the probability of having two sample points in any given straight line is zero. 

While $S_n={\mathcal X}_n$ is indeed a very simple estimator of $S$ it is also obviously unsatisfactory in many important aspects. In particular, except for trivial situations, it will typically fail to converge with respect to the distance ``in measure'' \eqref{dmu}. Also,   it will not give in general a consistent estimator of $\partial S$ since in  $d_H(\partial S_n,\partial S)\nrightarrow 0$ a.s.,  except for some particular distributions (e.g., discrete distributions with a bounded support).

\subsection{The $\theta$-biconvex hull}
 
Theorem \ref{th:p2} suggests a  natural way to get a meaningful non-trivial biconvex estimator of $S$ when the axes are known (let us denote them $b_1,b_2$). Indeed, since this theorem establishes that, whenever $\partial S$ is path connected, $S$ is biconvex if and only if $S$ fulfils a particular case of the $\pi/2$-lighthouse property.   This lead us to use the following version of the hull notion as an estimator of $S$ (recall that, given a unit vector $v$ and $y\in{\mathbb R}^2$, $C_v(y)$ denotes the open cone with vertex $y$, axis in the direction $v$  and opening angle $\pi/2$).

\begin{definition}\label{def:bh}
Given a set $A\subset{\mathbb R}^2$ and an angle $\theta\in[0,\pi/2)$, let $b_1=(\cos(\theta),\sin(\theta))$ and $b_2=(-\sin(\theta),\cos(\theta))$, where the coordinates are in the canonical basis, and  $\xi=(b_1+b_2)/\|b_1+b_2\|$. We define the lighthouse  $\theta$-biconvex hull (or just the $\theta$-biconvex hull) of $A$ by
	\begin{equation*}\label{lb2}
	{\mathbb B}_\theta(A)=\bigcap\big\{(C_{R^i(\xi)}(y))^c:\ i=0,1,2,3,\ y \mbox{ is such that } C_{R^i(\xi)}(y)\cap A=\emptyset\big\},
	\end{equation*}

\end{definition}

The intuitive idea behind Definition \ref{def:bh} is quite clear: let us consider all possible open quadrants (i.e. cones with opening angle $\pi/2$) whose sides are $\theta$-half lines or $\pi/2$-rotations of such half lines. Then, the $\theta$-biconvex hull of a set $A\subset{\mathbb R}^2$ is just the intersection of the complements of all quadrants of this type which do not intersect $A$.

\begin{remark}\label{remark:r4} 
 Note that, as a consequence of Theorem \ref{th:p2}, if S is $\theta$-biconvex, then ${\mathbb B}_\theta(S)=S$.
We will mainly use the notion of $\theta$-biconvex hull  for the particular case where $A={\mathcal X}_n$ is a finite sample of points randomly drawn on a $\theta$-biconvex set $S$. In that case
${\mathbb B}_\theta({\mathcal X}_n)$ is used as an estimator of $S$. We will show in Theorem \ref{thhull} that this is indeed a reasonable estimator for $S$, under some regularity conditions. When the angle $\theta$ is given and fixed we will sometimes omit the sub-index $\theta$ in ${\mathbb B}_\theta({\mathcal X}_n)$.

The notion of biconvex hull introduced in Definition \ref{def:bh} has been already considered (with a completely different motivation) in \citet[Def. 2.4]{ott84} under the name of ``maximal rectilinear convex hull''. These authors also outline an algorithm to evaluate ${\mathbb B}({\mathcal X}_n)$.
However, for our statistical purposes we will propose another slightly different algorithm to construct ${\mathbb B}({\mathcal X}_n)$. In fact, the idea behind this algorithm will be used later (see Theorem \ref{thhull}) to prove some relevant statistical properties of ${\mathbb B}({\mathcal X}_n)$ as an estimator of $S$. 

\end{remark}
 The following theorem establishes a natural property of the biconvex hull: the $\theta$-biconvex hull of a set $S$ which is not $\theta$-biconvex is strictly larger (in measure) than the original set. This property will be useful later in the proof of Theorem \ref{th:anglechoice}.

\begin{theorem} \label{th_4_3} Let $S$ be in the hypotheses of Theorem \ref{th:p2}.  If $S$ is $\theta_0$-biconvex but not $\theta$-biconvex for some $\theta_0, \theta\in [0,\pi/2)$ with $\theta_0\neq \theta$, then $\mu({\mathbb B}_\theta(S)\setminus S)>0$.
\end{theorem}

\begin{proof} Let us assume without loss of generality that $\theta_0=0$ (if this is not the case change the canonical axes).   Since $S$ is assumed to be $\theta_0$-biconvex, but not $\theta$-biconvex,  the equivalence between biconvexity and cone supporting property established in Theorem \ref{p1} holds for all points in $\partial S$ with the angle $\theta_0$ (and axis $\xi_0$) and fails for some $x\in \partial S$ with the angle $\theta$ (and axis $\xi$). In other words, there exists $x\in\partial S$ such that, for some $i\in\{0,1,2,3\}$, the quadrant $C_{R^i(\xi_0)}(x)$ fulfills $C_{R^i(\xi_0)}(x)\cap S=\emptyset$, but for all $j=0,1,2,3$, $C_{R^j(\xi)}(x)\cap S\neq\emptyset$.   We can assume without loss of generality that $i=0$.  Let us assume also that $j=1$, the other cases are treated similarly. Hence, there exists some $a\in S\cap C_{R(\xi)}(x)$ (since the $\theta$-biconvexity condition established in Theorem \ref{p1} fails at $x$) and also it does exist some $d\in S\cap C_{R^3(\xi_0)}(x)\cap C_{R^3(\xi)}(x)$; to see this note that, as we are assuming $\theta_0=0$, the set $\partial S\cap C_{R^3(\xi)}(x)$ (which is not empty) has only some overlapping with either $C_{R^3(\xi_0)}(x)$ or $C_{R^0(\xi_0)}(x)$. But, by assumption, $S\cap C_{R^0(\xi_0)}(x)=\emptyset$.
We can also assume that $a,d\neq x$ and they are in $\partial S$ since, for instance, if $a \in int(S)$, the line passing throughout $x$ and $a$ meets $\partial S$ at a point $a'$ such that $\|a'\|>\|a\|$, then define $a=a'$. Observe that this point does not belong to the boundary of the cone  $C_{R(\xi)}(x)$.
This proves that  $\ang (a-x,d-x)<\pi$. Still,  by construction, $\ang (a-x,d-x)\geq \pi/2$.  Denote ${\mathbf c}_1$, ${\mathbf c}_2$ the axes  passing through $x$  determined by $\theta$  (that is the $\theta$-line is the bisectrix of one of the quadrants determined by the vectors ${\mathbf c}_i$, $i=1,2$)  in such a way that, if $a=(a_1,a_2)$, $d=(d_1,d_2)$ stand for the coordinates of $a$ and $d$ wrt ${\mathbf c}_1$ and ${\mathbf c}_2$, we have $d_1>0,d_2>0$.   Let us assume that $a_2\geq d_2$. Since  $\ang (a-x,d-x)\in[\pi/2,\pi)$, we must have   $a_1\leq 0$ , $a_2>0$. 
	Let us denote $z,y$ the two intersection points of the line containing $d$, parallel to $\mathbf{c}_1$  (see Figure \ref{mudif}). 
	 We are going to prove that the open triangle, $E$, determined by $y,x,z$ is included in ${\mathbb B}_\theta(S)$. By construction $E\subset S^c$. Let $\gamma\subset \partial S$ a curve joining $a$ and $d$. Note that $\gamma\subset E^c$. Consider $t\in E$. Since $\gamma\subset E^c$ and $a_2\geq d_2$, $C_{R^{2}(\xi)}(t)\cap\gamma\neq \emptyset$ and $C_{R^3(\xi)}(t)\cap \gamma\neq \emptyset$. Then $t\in {\mathbb B}_\theta(S)$. Finally  $\mu({\mathbb B}_\theta(S)\setminus S)\geq \mu(E)>0$.

	\begin{figure}[!ht]
		\centering
		\includegraphics[scale=0.6]{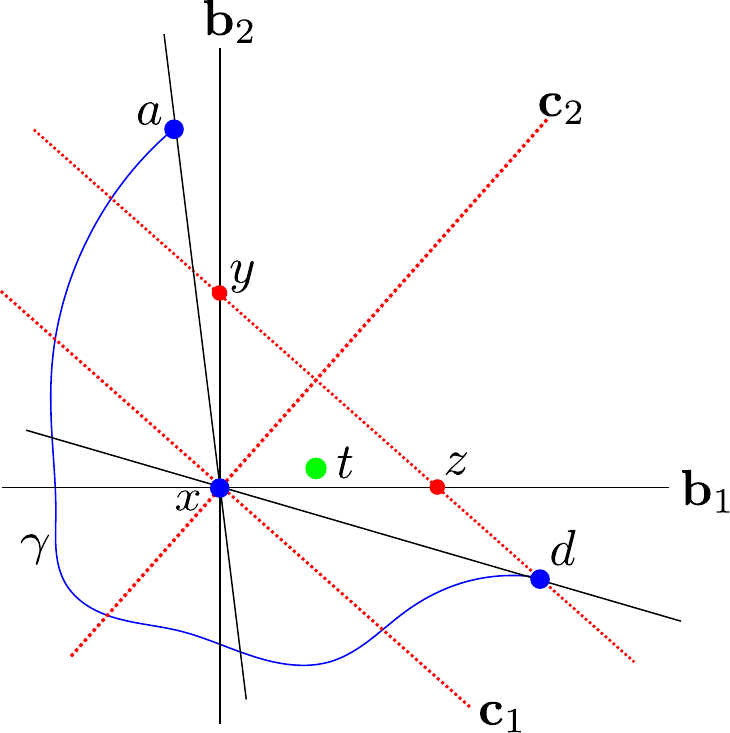} \ \ \ \ 
		\caption{The triangle determined by $a,x,d$ is included in ${\mathbb B}_\theta(S)$.}
		\label{mudif}
	\end{figure}
\end{proof}

The following theorem is the most important statistical result of this paper. It is concerned with the convergence rates properties of the sample biconvex hull as an estimator of a biconvex set.  Note, however, that the result is not ``statistical'' itself in the sense that it is a bit more general as it concerns the approximation of a biconvex set $S$ by a finite number of points (not necessarily random) inside $S$.

\begin{theorem} \label{setestth} Let $S$ in the hypotheses of Theorem \ref{th:p2}. Let $\X_n$ be any set of $n$ points in $S$ (not necessarily random). Then 
for all $\theta\in [0,\pi/2)$,  for all $n$ and for all $c>\sqrt{2}$

\begin{equation} \label{bordes2}
d_H({\mathbb B}_\theta(\X_n),{\mathbb B}_\theta(S))\leq cd_H(S,\X_n),
\end{equation} 
\begin{equation} \label{bordes}
d_H(\partial {\mathbb B}_\theta(S), \partial {\mathbb B}_\theta(\X_n))\leq cd_H(S,\X_n).
\end{equation}
 If for all $\theta$ there exists there the Minkowski content of $\partial {\mathbb B}_\theta(S)$, $L_0=L_0(\theta)$, then, for $n$ large enough, 
\begin{equation} \label{dnu}
d_{\mu}({\mathbb B}_\theta(\X_n),{\mathbb B}_\theta(S))\leq 2cL_0(\theta) d_H(S,\X_n) .
\end{equation} 
In particular, if, for every $n\in{\mathbb N}$, $\X_n$ denotes a random sample from a distribution $P_X$ with a $\theta$-biconvex support $S$, expressions \eqref{bordes}, \eqref{bordes2} and \eqref{dnu} provide almost sure consistency results for the estimation of $S$ (with respect to $d_H$ and $d_{\mu}$) and $\partial S$ (with respect to $d_H$).
\end{theorem}
\begin{proof}
    Let us prove that for all $n$, if we denote $\eps_n=d_H(S,\X_n)$, for all $\theta \in [0,\pi/2)$ and for all $\eta>0$
		\begin{equation} \label{inc1}
		{\mathbb B}_\theta(S)\ominus  B(0,(\sqrt{2}+\eta)\eps_n)\subset {\mathbb B}_\theta(\X_n)\subset {\mathbb B}_\theta(S).
		\end{equation}
		Let us prove the first inclusion. If  ${\mathbb B}_\theta(S)\ominus  B(0,(\sqrt{2}+\eta)\eps_n)=\emptyset$ then \eqref{inc1} holds trivially. Otherwise,   by contradiction   let $x\in {\mathbb B}_\theta(S)\ominus B(0,(\sqrt{2}+\eta)\eps_n)$ but $x\notin {\mathbb B}_\theta(\X_n)$. Since $x\notin {\mathbb B}_\theta(\X_n)$  we have, from Definition \ref{def:bh}, that  there exists $\xi=\xi(\theta)$ with $\|\xi\|=1$ and $C_{R^i(\xi)}(y)$ for some $i=0,1,2,3$, such that  $x\in C_{R^i(\xi)}(y)$ and $C_{R^i(\xi)}(y)$ is disjoint with $\X_n$. Let us assume,  without loss of generality, that this holds for $i=0$, so that $x\in C_{\xi}(y)$, $C_{\xi}(y)\cap \X_n =\emptyset$. Since $x\in {\mathbb B}_\theta(S)\ominus B(0,(\sqrt{2}+\eta)\eps_n)$ we have that $B(x,(\sqrt{2}+\eta)\eps_n)\subset {\mathbb B}_\theta(S)$. Observe that $\sin(\pi/4)=1/\sqrt{2}$, then there exists $z\in B(x,(\sqrt{2}+\eta)\eps_n)\cap C_{\xi}(y)$ such that $d\big(z,\partial C_{\xi}(y)\big)> \eps_n$ (see Figure \ref{th3_2}). Since $z\in {\mathbb B}_\theta(S)$ 
		we have that $\overline{C_{\xi}(z)}\cap S\neq \emptyset$. Let $s\in \overline{C_{\xi}(z)}\cap S$, $\overline{B(s,\eps_n)}\subset C_{\xi}(y)$, since $\eps_n=d_H(S,\X_n)$ this contradicts  $C_{\xi}(y)\cap \X_n=\emptyset$. From (\ref{inc1}) it follows (\ref{bordes2}).\\
	  To prove (\ref{bordes}), let us assume by contradiction that there exists $x\in \partial {\mathbb B}_\theta(\X_n)$ and $x\in {\mathbb B}_\theta(S)\ominus B(0,(\sqrt{2}+\eta)\eps_n)$. Since $x\in \partial {\mathbb B}_\theta(\X_n)$,   from Theorem \ref{p1},  there exists $\xi=\xi(\theta)$ with $\|\xi\|=1$ such that $C_{R^i(\xi)}(x)\cap \X_n=\emptyset$ for some $i$. To simplify the notation assume $i=0$. By \eqref{inc1} for all $n$ and all $\eta>0$, $B(x,(\sqrt{2}+\eta)\eps_n)\subset {\mathbb B}_\theta(S)$. Let us consider the cones $C_{\xi}(x+\lambda\xi) \quad 0\leq\lambda<(\sqrt{2}+\eta)\eps_n,$
		then $C_{\xi}(x+\lambda\xi) \cap S\neq \emptyset$ for all $0\leq\lambda<(\sqrt{2}+\eta)\eps_n$. Since $S$ is compact, $\overline{C_{\xi}(x+(\sqrt{2}+\eta)\eps_n\xi)} \cap S\neq\emptyset$. Let $s\in\overline{C_{\xi}(x+(\sqrt{2}+\eta)\eps_n\xi)}\cap S$. By the geometric argument we made above, it follows that the distance from $z$ to the boundary of the ``inner cone'' $C_{\xi}(x+(\sqrt{2}+\eta)\eps_n\xi)$ is at least $\epsilon_n$; thus the distance from $z$ to  $\partial C_{\xi}(x)$ is at least $2\eps_n$ and we have $B(s,2\eps_n)\subset C_{\xi}(x)$ 
		but, since $C_{\xi}(x)\cap \X_n=\emptyset$,  this a contradiction with $\epsilon_n=d_H(S,\X_n)$.\\

		Let us prove \eqref{dnu}. Observe that, from \eqref{inc1}, for all $n>0$, for all $\theta$ and for all $c>\sqrt{2}$,
		$$d_{\mu}({\mathbb B}_\theta(\X_n),{\mathbb B}_\theta(S))\leq \mu\big(B(\partial {\mathbb B}_\theta(S),c\varepsilon_n)\big),$$
		where $\varepsilon_n=d_H(S,\X_n)$.  
		From  $L_0(\theta)=\lim_{\epsilon\to 0} \mu(B(\partial {\mathbb B}_\theta(S),\epsilon)/(2\epsilon)$, we get that, for $n$ large enough
		$$
		\mu(B(\partial {\mathbb B}_\theta(S),c\epsilon_n))\leq 
		2c'L_0(\theta) d_H(S,\X_n)
		$$
		for all $c'>c$, which in  turn implies  \eqref{dnu}. The final claim in the statement follows directly since, from Borel-Cantelli Lemma $d_H(\X_n,S)\to 0$, a.s. 
	 	\begin{figure}[!ht]
	 	\centering
	 	\includegraphics[scale=0.8]{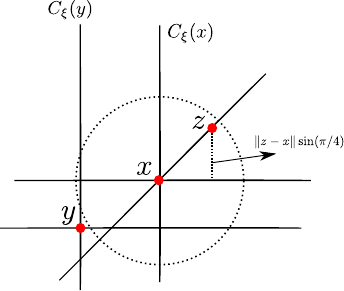} \ \ \ \ 
	 	\caption{}
	 	\label{th3_2}
	 \end{figure}
		
\end{proof}

\

\noindent \textit{Convergence rates }

\

 As a consequence of Theorem \ref{setestth}, we can easily derive convergence rates, under an additional shape condition of ``standardness'' for the set $S$. This shape restriction is quite popular 
in set estimation; see \cite{cue97}, \cite{rin10}. The formal definition is as follows. 
\begin{definition} \label{standard} A set $S\subset{\mathbb R}^2$ is said to be \it standard\/ \rm with respect to a Borel measure $\nu$ if there exist $\lambda>0$, $\delta>0$ such that
	\begin{equation}\label{estandar}
	\nu(B(x,\epsilon)\cap S)\geq \delta\mu(B(x,\epsilon)),\ \mbox{for all } x\in S, 0<\epsilon\leq \lambda.
	\end{equation}
\end{definition}

\begin{corollary}\label{cor:rates}
  Under the hypotheses of Theorem \ref{setestth}, if we additionally assume that $S$ is standard with respect to $P_X$, we have that the left-hand sides distances in expressions \eqref{bordes}, \eqref{bordes2} and \eqref{dnu} tend to 0, at a rate of type  $\mathcal{O}(\log(n)/n)^{1/d}$ (a.s.) as $n$ tends to infinity.  
\end{corollary}
\begin{proof}
The result follows from the fact that, if $S$ is compact and standard with respect to $P_X$, $d_H(S,\X_n)=\mathcal{O}(\log(n)/n)^{1/d}$ (see Theorem 3 in \cite{cue04}). 
\end{proof}

 It is worth noting that the convergence rate $d_H(S,\X_n)=\mathcal{O}(\log(n)/n)^{1/d}$ is the same rate obtained for the estimation of the convex hull as well as the $\rho$-cone convex hull (see \cite{chola:14}) for details.

\subsection{An algorithm to construct the sample $\theta$-biconvex hull }\label{knax}

 Throughout this section we assume that ${\mathcal X}_n=\{X_1,\ldots,X_n\}$  is a sample from a distribution $P_X$ with compact support $S\subset{\mathbb R}^2$. 
We will provide below an exact algorithm to build the $\theta$-biconvex hull, ${\mathbb B}_\theta({\mathcal X}_n)$ 
with edges along the directions $b_1,b_2$, given by $b_1=(\cos(\theta),\sin(\theta))$, $b_2=(-\sin(\theta),\cos(\theta))$). Such algorithm directly relies upon Definition \ref{def:bh}. Thus, our goal is the intersection of the complements of all (open) $\pi/2$-cones of type $C_{R^i(\xi)}(x)$, disjoint with ${\mathcal X}_n$, whose axes correspond to either the direction $\xi=(b_1+b_2)/\|b_1+b_2\|$ or any of its $\pi/2$-rotations, $R^i(\xi)$ ($i=0,\ldots,4$).

The basic idea is as follows. The biconvex hull must be necessarily contained in the ordinary convex hull $E=co({\mathcal X}_n)$. Thus, we take this set as a starting element. Then, for each sample point $X_i$, let us select among the four (open) cones $C_i\equiv C_{R^j(\xi)}(X_i)$ with vertex $X_i$ (and axes with the directions $R^j(\xi)$) those ``empty cones'' $C_i$ not including any other sample point. For every such ``empty'' cone, let us consider the ``maximal empty horizontal cone'' $C_i^1$, obtained by horizontally moving $C_i$ until some other sample point is met. Similarly, calculate the ``maximal vertical cone'' $C_i^2$ obtained by vertically moving $C_i$ until we met some other sample point. Then, take out from $E$ the intersections $E\cap C_i^1$ and $E\cap C_i^2$. Iterate this process for all the remaining sample points. The biconvex hull of the sample ${\mathcal X}_n$ is what is left of $E$ after such deletion process. 

In more schematic terms, the algorithm is as follows.

\begin{itemize}
    \item[START:] Put $E=co({\mathcal X}_n)$.
    \item [ITERATION:] For each $i=1,\ldots,n-1$ calculate the four cones $C_{R^j(\xi)}(X_i)$, $j=0,\ldots,3$. If all these cones contain sample points in ${\mathcal X}_n$, then put $i\leftarrow i+1$ and repeat the process of cones calculation. Otherwise,
    \begin{itemize}
        \item[I1.] For each ``empty'' cone $C_i$ (i.e., a cone not containing sample points) calculate the ``horizontal maximal empty cone'' $C_i^1$, obtained by moving horizontally $C_i$ until some other sample point is met. Replace $E\leftarrow E\setminus C_i^1$. 
        \item[I2.] Calculate also the ``vertical maximal empty cone'' $C_i^2$, obtained by moving vertically $C_i$ until some other sample point is met. Replace $E\leftarrow E\setminus C_i^2$.
    \end{itemize} 
    \item[OUTPUT:] The  set $E:=E_n$ obtained from the above process after all points have been considered in the iterations. 
\end{itemize}

\

Observe that in every step of the algorithm we remove cones (that is, we intersect with the complement of a cone, which is a biconvex set), and the intersection of biconvex sets is also biconvex, then the result of the algorithm is a biconvex set.  In Theorem \ref{thhull} below we will show that, in fact the algorithm output coincides with ${\mathbb B}_\theta(\X_n)$, the $\theta$-biconvex hull of the sample. 
Other relevant statistical properties of the set $E_n={\mathbb B}_\theta(\X_n)$ when considered as an estimator of an unknown $\theta$-biconvex compact set $S$ are also established. Of course, this makes sense in the case that $\X_n$ is a random sample drawn from a probability distribution with support $S$. 

\

\begin{theorem} \label{thhull} Let $S$ in the hypotheses of Theorem \ref{th:p2}. Let $\X_n$ be any set of $n$ points in $S$.   For all $\theta\in [0,\pi/2)$,  the final output of the above algorithm coincides with the biconvex hull, that is, $E_n={\mathbb B}_\theta(\X_n)$.

\begin{proof} 
By definition, ${\mathbb B}_\theta(\X_n)\subset E_n$. To prove the other inclusion, let $x\notin {\mathbb B}_\theta(\X_n)$ we will prove that $x\notin E_n$. Since $x\notin {\mathbb B}_\theta(\X_n)$ there exist a cone $C_{R^i(\xi)}(y)$ for some $i=0,1,2,3$, with $x\in C_{R^i(\xi)}(y)$ and $C_{R^i(\xi)}(y)\cap \X_n=\emptyset$.  Denote by $l_1$, $l_2$ the two half-lines defining the boundary of $C_{R^i(\xi)}(y)$. If there are points of $\X_n$ in both $l_1$ and $l_2$ then, by construction of $E_n$,  we have $x\in E_n^c$. Otherwise,  translate the cone until meeting the sample in both half-lines of the translated boundary. Then the result follows by applying again the above argument to the translated cone. 
\end{proof}
\end{theorem}

\

 The following continuity  result has some conceptual and practical interest. 

\begin{corollary} \label{cr_revisado} Let $S$ be in the hypotheses of Theorem \ref{th:p2}. Assume further that
	\begin{itemize}
		\item[(a)] For all $\theta\in[0,\pi/2]$ the limit $$
		L_0(\theta)=\lim_{\epsilon\rightarrow 0}\frac{\mu\left(B(\partial{\mathbb B}_\theta(S),\epsilon)\right)}{2\epsilon}
		$$
		is finite and uniform on $\theta$.
		\item[(b)] There exists $L>0$ such that $L_0(\theta)\leq L$ for all $\theta$.		
	\end{itemize}
	Then, the function $\theta\mapsto \mu( {\mathbb B}_\theta(S))$ is continuous. 
	
	\begin{proof} 
For each $n$, take a set denoted by $\mathcal{X}_n$ of $n$ points included in $S$ in such a way that $d_H(S,\X_n)\to 0$. Recall that in Theorem \ref{setestth}, equation \eqref{dnu}, we proved that, given $\epsilon>0$ there exists,  $n_0=n_0(\theta)$ such that, for all $c>\sqrt{2}$,
		\begin{equation}
		d_{\mu}({\mathbb B}_\theta(\X_n),{\mathbb B}_\theta(S))\leq2cL_0(\theta) d_H(S,\X_n), \quad \mbox{for all } n>n_0. \label{eq:n0}
		\end{equation}
		If we revise the proof of this inequality, we can readily see that, under assumptions (a) and (b) above, a inequality of type 
		\begin{equation}
		d_{\mu}({\mathbb B}_\theta(\X_n),{\mathbb B}_\theta(S))\leq 2cL d_H(S,\X_n), \quad \mbox{for all } n>n_1. \label{eq:n1}
		\end{equation}
		holds for some index $n_1$ not depending on $\theta$. In other words, inequality \eqref{dnu}, where $L_0(\theta)$ is replaced with $L$,  holds uniformly on $\theta$. Hence
		\begin{align}
		&\sup_\theta|\mu({\mathbb B}_\theta(\X_n))-\mu({\mathbb B}_\theta(S))|\leq \sup_\theta d_{\mu}({\mathbb B}_\theta(\X_n),{\mathbb B}_\theta(S)) \\ \nonumber
		&\leq 2cL d_H(S,\X_n)\rightarrow 0,\label{eq:n2}
		\end{align}
	by construction of $\X_n$. Now, note that the transformations $\theta\mapsto \mu({\mathbb B}_\theta({\mathcal X}_n))$ are continuous. To see this, recall that, according to Theorem \ref{thhull},  $E=\mathbb{B}_\theta(\X_n)$ where $E$ is constructed as the intersection of the complements of a finite number of quadrants meeting some sample points at their boundaries. Then, by construction, $\theta\mapsto\mu({\mathbb B}_\theta({\mathcal X}_n))$ is a continuous function of $\theta$. 
		
		As a conclusion, $\theta\mapsto \mu({\mathbb B}_\theta(S))$ is also continuous, as it can be expressed   as a uniform limit of continuous functions $\theta\mapsto\mu({\mathbb B_\theta}({\mathcal X}_n))$.

	\end{proof}
\end{corollary}	

\begin{remark}
	Regarding assumption (a) in Corollary \ref{cr_revisado} note that it is automatically fulfilled, whenever all the sets $\partial {\mathbb B}_\theta(S)$ have a linear volume function   with  bounded coefficients,  that is, there exist a bounded function $L_0(\theta)$ and a constant $R>0$ such that
	\begin{equation}\label{eq:pv}
	\mu(B(\partial {\mathbb B}_\theta(S),\epsilon))=2L_0(\theta)\epsilon,\ \mbox{for all } \epsilon\in [0,R].
	\end{equation}
	In particular, an expression of type \eqref{eq:pv} holds whenever $\partial {\mathbb B}_\theta(S)$ has a polynomial volume and is homeomorphic to the unit circle $\partial B(0,1)$.  See \cite{cue18},  and references therein, for a detailed account of the meaning of the polynomial volume assumption and its statistical applications. 

\end{remark}

\subsection{When the biconvexity axes are unknown}

 If $S$ is  $\theta$-biconvex but the value of $\theta$ is unknown, we can still approximate $\theta$ from a data-driven sequence $\hat\theta_n$. The definition of this approximating sequence and the sense in which it approaches the true $\theta$   is made explicit in the following theorem.

\begin{theorem}\label{th:anglechoice}
Under the hypotheses on $S$ imposed in Corollary \ref{cr_revisado}, denote by $A$ the set  of angles $\theta$ in $[0,\pi/2]$ for which the set $S$ is $\theta$-biconvex. Assume that $A$ is non-empty. For $\theta\in [0,\pi/2]$ define the sequence of random functions $\Psi_n(\theta)=\mu({\mathbb B}_{\theta}(\X_n))$. 
Let $\{\hat \theta_n\}$ be 
a sequence of random variables such that
\begin{equation}\label{argmin1}
\hat{\theta}_n\in \mbox{argmin}_\theta\Psi_n(\theta)
\end{equation}
Then, with probability one all the accumulation points of the sequence $\{\hat \theta_n\}$  of minimizers of $\Psi_n$,  belong to $A$.
	\begin{proof}
		First note that, by construction, ${\mathbb B}_{0}(\X_n)={\mathbb B}_{\pi/2}(\X_n)$. Note also that the set $A$ of biconvexity angles is compact; this follows directly from the continuity of the function $g(\theta)=\mu({\mathbb B}_{\theta}(S))-\mu(S):=\Psi(\theta)-\mu(S)$ (see Corollary \ref{cr_revisado}),   together with $A=g^{-1}(0)$ (this follows from  Remark \ref{remark:r4}  and Theorem \ref{th_4_3}).  Thus, $A$ is a closed set included in the compact set $[0,\pi/2]$ and therefore compact. Now, reasoning by contradiction, suppose that, with positive probability, there is a subsequence  of a sequence $\{\hat \theta_n\}$ of minimizers, denoted again $\{\hat \theta_n\}$ for simplicity, converging to a point  $\theta_1\in A^c$. From Theorems \ref{th:p2} and \ref{setestth}, if $S$ is $\theta$-biconvex, that is, if $\theta\in A$, then ${\mathbb B}_\theta(S)=S$ and $\Psi_n(\theta)\rightarrow \mu(S)$ a.s., while if $S$ is not $\theta_1$-biconvex,  from Theorems \ref{th_4_3} and \ref{setestth}, $\mu({\mathbb B}_{\theta_1}(\X_n))-\mu(S)\rightarrow \mu({\mathbb B}_{\theta_1}(S)\setminus S)>0$ a.s. In any case,  from  the proof of Corollary \ref{cr_revisado}, we know that $\Psi_n(\theta)\to \Psi(\theta)$ a.s. uniformly on $\theta$. 
		
		Now, take $\epsilon>0$ small enough so that there is a closed neighbourhood $U(\theta_1)$ of $\theta_1$ in $[0,\pi/2]$ such that $U(\theta_1)\subset A^c$ and (using the continuity of $\Psi$),  $\Psi(\theta)>\mu(S)+\epsilon$ for all $\theta\in U(\theta_1)$. 
		
		Also, since  $\Psi_n(\theta)\to \mu(S)$, uniformly on $\theta\in A$, a.s., we have, for all $\theta\in A$,
		\begin{equation}\label{eq:theta211}
		\Psi_n(\hat \theta_n)\leq\Psi_n(\theta)<\Psi_n(\theta_1)-\epsilon/2\mbox{ eventually, almost surely}. 
		\end{equation}	
		But, on the other hand, the fact $\hat \theta_n\to \theta_1$ with positive probability and the a.s. uniform convergence $\Psi_n\to \Psi$ on $U(\theta_1)$ entail $\Psi_n(\hat\theta_n)\to \Psi(\theta_1)$ with positive probability, which contradicts \eqref{eq:theta211}.
		
	\end{proof}
\end{theorem}
 The following result concerns the estimation of a biconvex set with an estimated biconvexity direction.

\

\begin{theorem}\label{th:consistency1} Let $S$ be in the hypotheses of Theorem \ref{th:anglechoice}, let  $\hat{\theta}_{n}$ be any  sequence of angles fulfilling  \eqref{argmin1}. Then, we have
	\begin{equation} \label{convdmu1}
	d_{\mu}({\mathbb B}_{\hat{\theta}_{n}}(\X_{n}),S)\to 0 \quad a.s.
	\end{equation}

\end{theorem}
\begin{proof}   
	First note that $\hat \theta_n$ is a random variable (defined on some probability space $(\Omega,{\mathcal A},{\mathbb P})$) depending on the sample points $X(\omega),\ldots,X_n(\omega)$. We will show that \eqref{convdmu1} holds for almost every $\omega\in\Omega$. Then, take any $\omega$ such that the convergence in measure from $\X_n(\omega)$ to $S$ established in \eqref{eq:n1} holds valid and put $\X_n(\omega)=\X_n$ and $\hat \theta_n=\theta_n(\omega)$. 
	
	Now, note that in order to establish \eqref{convdmu1}, it suffices to show that any subsequence $\{d_{\mu}({\mathbb B}_{\hat{\theta}_{n_k}}(\X_{n_k}),S)\}\subset \{d_{\mu}({\mathbb B}_{\hat{\theta}_{n}}(\X_{n}),S)\}$ contains a further subsequence converging to 0.
	
	As $\{\hat \theta_{n_k}\}$ is a 
	sequence included in the compact set $[0,\pi/2]$. there is a further subsequence (denoted again $\{\hat \theta_{n_k}\}$ by simplicity) 
	convergent to some $\theta_0$. 
	From Theorem \ref{th:anglechoice} we must have $\theta_0\in A$. Then, 
	we have
	\begin{align}\label{eq:triangular}
	&d_{\mu}({\mathbb B}_{{\hat \theta}_{n_k}}(\X_{n}),S)=d_{\mu}({\mathbb B}_{\hat{\theta}_{n_k}}(\X_{n_k}),{\mathbb B}_{\theta_0}(S))\\ \nonumber
	&\leq d_{\mu}({\mathbb B}_{\hat{\theta}_{n_k}}(\X_{n_k}),{\mathbb B}_{\hat \theta_{n_k}}(S))+d_{\mu}({\mathbb B}_{\hat \theta_{n_k}}(S),{\mathbb B}_{\theta_0}(S))
	\end{align}
	Now, given $\epsilon>0$, the term $d_{\mu}({\mathbb B}_{\hat{\theta}_{n_k}}(\X_{n_k}),{\mathbb B}_{\hat \theta_{n_k}}(S))$ can be made smaller than $\epsilon/2$ for $n$ large enough. This follows from expression \eqref{eq:n1} in the proof of Corollary \ref{cr_revisado}. The second term $d_{\mu}({\mathbb B}_{\hat \theta_{n_k}}(S),{\mathbb B}_{\theta_0}(S))$ is also eventually smaller than $\epsilon/2$, as a consequence of the uniform continuity of $\theta\mapsto \mu({\mathbb B}_\theta)$, since $$d_{\mu}({\mathbb B}_{\hat \theta_n}(S),{\mathbb B}_{\theta_0}(S))=\mu({\mathbb B}_{\hat \theta_n}(S)\setminus S)=\mu({\mathbb B}_{\hat \theta_n}(S))-\mu({\mathbb B}_{ \theta_0}(S)).$$
	Thus, we have proved that for any subsequence extracted from  $\{d_{\mu}({\mathbb B}_{\hat{\theta}_{n}}(\X_{n}),S)\}$ there is a further subsequence converging to 0 and the proof is complete.

\end{proof}

\section{ Some examples}\label{sec:sim}
A few final examples   are included here just for illustrative purposes, in order to gain some intuition on the practical meaning of the notions we have introduced.
 
\subsection{Estimation of biconvex but not $\alpha$-convex set} 
Let us consider the set $S=T_1\cup T_2$,   where $T_1$ and $T_2$ are the closed triangles with vertices $\{(0,0), (1/2,1/2),(0,1)\}$ and $\{(1,0), (1/2,1/2),(1,1)\}$, respectively.   
Observe that $S$ is $\theta$-biconvex for $\theta=\pi/4$. 

We draw a uniform sample of $n$ data points over $S$ and we aim at reconstructing $S$ from such sample, using the $\pi/4$-biconvex hull of the data points (see Definition \eqref{lb2}). For comparison purposes, we also consider another usual set estimator, namely the $\alpha$-convex hull,  associated with the idea of $\alpha$-convexity above mentioned. Recall that the $\alpha$-convex hull of a sample of points is the intersection of the complements of all balls of radius $\alpha$ not including any sample point. See   \cite{pat10} for a description of the package \textit{alphahull} (which will be used in the example below to calculate $\alpha$-convex hulls); see  also \cite{cue12} for theoretical aspects and additional references on $\alpha$-convexity.

\begin{figure}[!ht]
	\centering
	\includegraphics[scale=0.32]{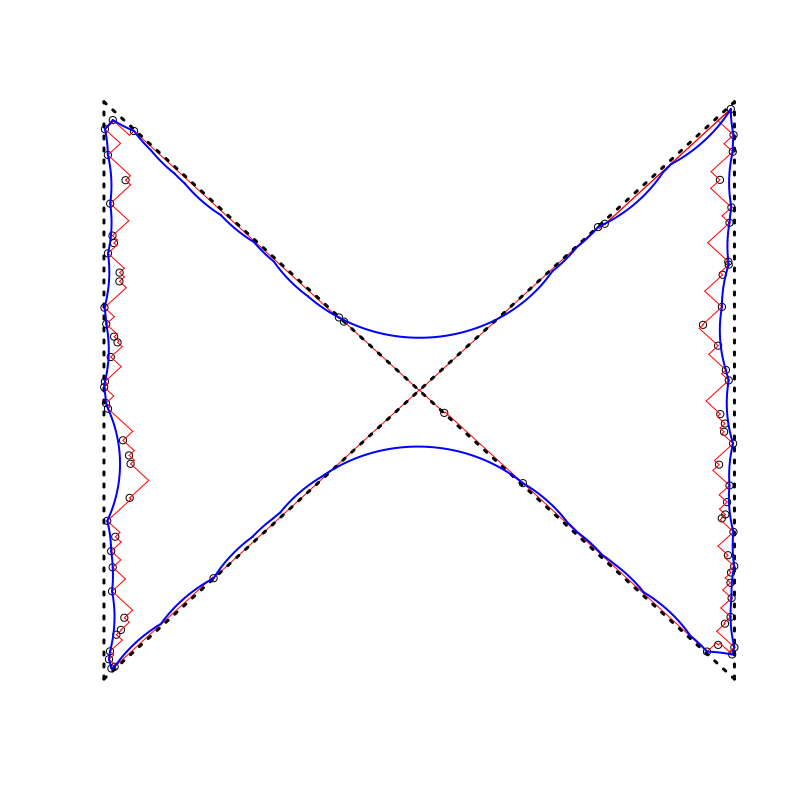}  \includegraphics[scale=0.32]{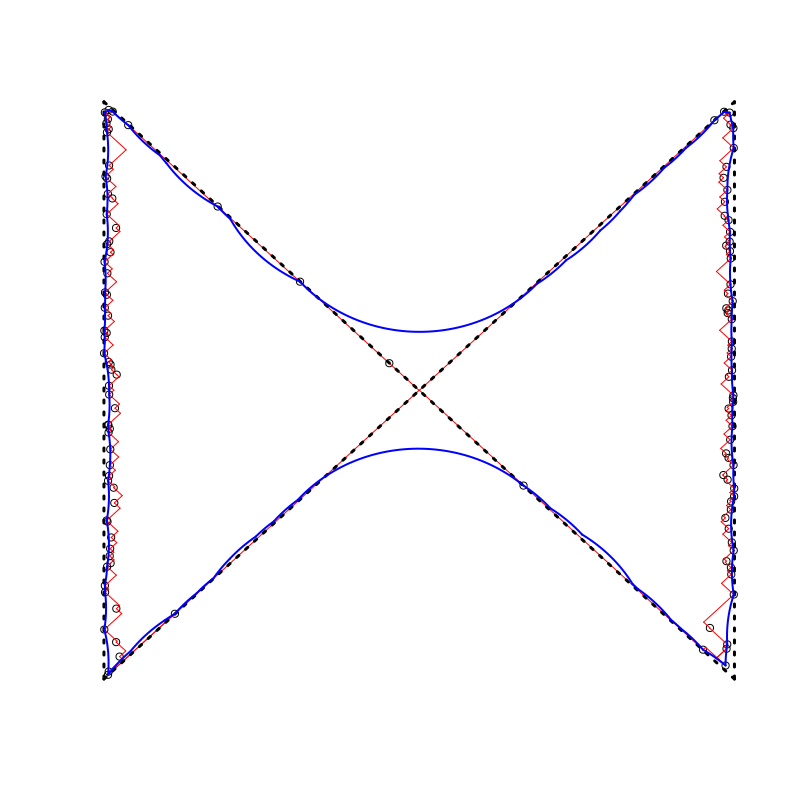}
	\caption{Estimation of a biconvex set $S$ using the biconvex hull and the $\alpha$-convex hull, with sample sizes $n=1000$ (left), and  $2000$ (right). The ``smoother'' line (in blue) is the boundary of the $1/3$ convex hull of the sample. The dotted line is $\partial S$. The more irregular line (in red) is the boundary of the biconvex hull of the sample.}
	\label{simu2}
\end{figure}

Figure \ref{simu2} shows both estimators of $S$ for the cases $n=1000$ (left) and $n=2000$ (right). 
The sample points are uniformly chosen over $S$.
The sample biconvex hull is  appears in the figure as  the set inside $S$ with a piece-wise linear boundary.  The sample $\alpha$-convex hull, with $\alpha=1/3$  is the set whose boundary is made of arcs of circles with radius $\alpha=1/3$.  

In this example the distance in measure (as defined in \eqref{dmu}, for $\nu=\mu$, the Lebesgue measure)  between the biconvex hull estimator and the true set was 0.04335, for $n=1000$. The analogous error measure for the $\alpha$-convex estimator
with $\alpha=1/3$ was 0.06746.

The respective values for $n=2000$ were 0.0319 and 0.06106. These distances have been approximated using a Monte Carlo sample of 50000 points drawn on $[0,1]^2$. We have also included in the comparison the $\epsilon_n$-offset of the sample points, defined by $\cup_{i=1}^n B(X_i,\epsilon_n)$. This is an all-purposes set estimator, sometimes called the Devroye-Wise (DW) estimator, which does not incorporate any prior shape information on $S$. It depends on a tuning parameter $\epsilon_n$. In this case, 
we have chosen $\epsilon_n=2.4(\log(n)/n)^{1/2}$, as suggested by Theorem 4 in \cite{cue04}. The results in Tables \ref{tab1} and \ref{tab2} show that, not surprisingly, the price to be paid for the generality of the DW estimator is some loss of efficiency when compared with the more specific estimators that incorporate some convexity-related information.
\begin{table}[h]

\begin{tabular}{cccc}
	$n$  &  $\mathbb{B}_{\pi/4}(\mathcal{X}_n)$ & $C_{1/3}(\mathcal{X}_n)$ &   DW   \\ \hline
	500  &                0.0686                &          0.0758          &  1.2015 \\
	1000 &                0.0444                &          0.0649          &  0.9962  \\
	1500 &                0.0392                &          0.0632          &  0.8431  \\
	2000 &                0.0275                &          0.0581          &  0.7361 \\
	2500 &                0.0273                &          0.0608          &  0.6667\\
\end{tabular}
\caption{Distance in measure averaged over 500 replications, for different values of $n$, for three set estimators.}
	\label{tab1}
\end{table}

 Of course, this example is, in some sense, ``favorable'' to the biconvex-hull estimator, since $S$ is a $\pi/4$-biconvex set but it is not $\alpha$-convex, for any $\alpha>0$ (since $S$ cannot possibly be expressed as the intersection of the complements of a family of open balls of any given positive radius). Thus our first small experiment should be seen as an assessment of how much improvement can be obtained in the estimation of $S$ by incorporating some shape information on $S$, not included in other better-known set estimators. 

\subsection{Estimation of a biconvex and $\alpha$-convex set}

\begin{figure}[!ht]
	\centering
	\includegraphics[scale=0.6]{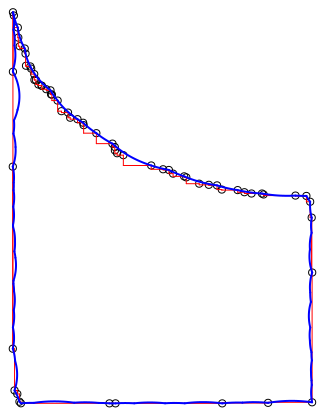} \hspace{2cm}
		\includegraphics[scale=0.57]{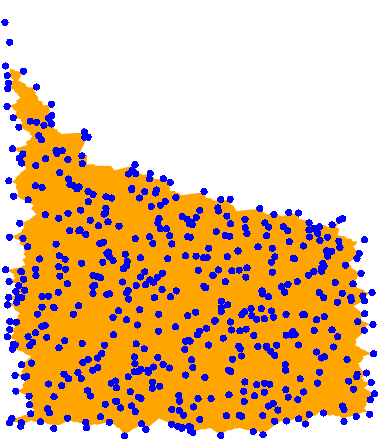}
	\caption{ Left panel: the ``wigglier'' line in red is the boundary of the biconvex hull (for $\theta=0$). The ``smoother'' line in blue is the boundary of the $1/3$-convex hull for a sample of size $n=2500$. Right panel: the orange-shaded area corresponds to the  $\pi/2$-cone convex hull by complements of the sample (with $n=500$)}
	\label{fig2}
\end{figure}

Our second example is based on the set  $S=[0,1]\times [0,2]\setminus B((1,2),1)$.  This set is  $\theta$-biconvex for several values of $\theta$, including $\theta=0$: in fact, we have chosen $\theta=0$ to construct biconvex hull of the sample. We have compared the biconvex hull of the sample with the $1/3$-convex hull, the 1-convex hull, the DW estimator (with the parameter $\epsilon_n$ chosen as before)  and the so-called $\pi/2$-cone-convex  hull by complements of the sample (see Definition \ref{lightcomp}), denoted by $\mathbb{C}_{\pi/2}(\mathcal{X}_n)$, as studied in \cite{chola:14}; see Figure \ref{fig2} right. This latter estimator must be calculated with an approximate stochastic algorithm, as described in \cite{chola:14} which has been constructed by removing 500 randomly selected cones.  The average over 500 replicates of the distance in measure between $S$ and the different estimators in competition is shown in Table \ref{tab2}, where again it can be seen that the all-purposes DW estimator is less efficient than those estimators incorporating shape information on $S$. 

\begin{table}[h]

\begin{tabular}{ccclcl}
	$n$  & $\mathbb{B}_{\pi/4}(\mathcal{X}_n)$ & $C_{1/3}(\mathcal{X}_n)$  &$C_{1}(\mathcal{X}_n)$  &   DW     & $\mathbb{C}_{\pi/2}(\mathcal{X}_n)$\\ \hline
	500  &                0.1181                &          0.1201          & 0.0782					&  1.2425  &  0.2010 \\
	1000 &                0.0530                &          0.0738          & 0.0424					&  1.1129  & 0.0148 \\
	1500 &                0.0431                &          0.0449          & 0.0310					&  0.9458  & 0.1228 \\
	2000 &                0.0343                &          0.0391          & 0.0305					&  0.8202  & 0.0987 \\
	2500 &                0.0353                &          0.0398          & 0.0233				    &  0.7316  & 0.0748 \\
\end{tabular}
\caption{Distance in measure averaged over 500 replications, for different values of $n$.}
\label{tab2}
\end{table}

In  Figure \ref{fig2} (left panel) we show, for $n=2500$, the boundary of the $1/3$-convex hull (the smoother line in blue) together with the boundary for $\mathbb{B}_{\pi/4}(\mathcal{X}_n)$   (the wigglier line in red) are shown.  The sample points drawn are those in the boundary of $\mathbb{B}_{\pi/4}(\mathcal{X}_n)$.  The right panel of Figure \ref{fig2} shows the sample points of a smaller sample (with $n=500$) together with the $\pi/2$-cone convex hull by complements (see Definition \ref{lightcomp}) of the sample, represented as the shaded area.

\subsection{Estimation of the biconvexity angle}
\begin{figure}[!ht]
	\centering
		\includegraphics[scale=0.45]{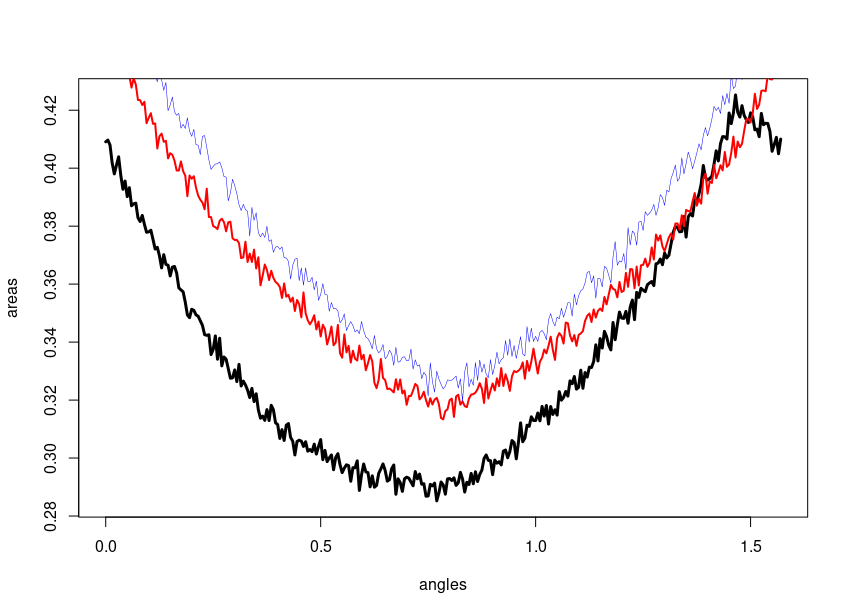}
	\caption{Representation of the curves $\Psi_n(\theta)= \mu(\mathbb{B}_\theta(\mathcal{X}_n))$ for $n=200, 400, 600$. The lower (resp. upper) curve corresponds to $n=200$ (resp. $n=600$).}
	\label{areas}
\end{figure}

Figure \ref{areas} shows the graphs of the functions $\Psi_n(\theta)= \mu(\mathbb{B}_\theta(\mathcal{X}_n))$ obtained for three sample sizes  ($n=200,400,600$). The value of $\theta$ varies on a grid from 0 to $\pi/2$ with 0.005 steps. The sample is uniformly distributed on the set $S=R_{\pi/4}([0,1]^2\setminus T)$, where $T$ is the triangle with vertices $(0,1)$, $(1/2,1/2)$ and $(1,1)$ and $R_{\pi/4}$ the clockwise rotation of angle $\pi/4$. Observe that this set is $\theta$-biconvex for $\theta=\pi/4\approx 0.78$. As a consequence of Theorems \ref{th_4_3} and \ref{th:anglechoice}, the convexity angle $\theta$ can be estimated by minimizing  $\Psi_n$. The lower curve corresponds to $n=200$, the intermediate one to $n=400$ and the upper one to $n=600$. The respective minima are attained at 0.77, 0.785 and 0.83.

\newpage

\begin{center}
    \sc Acknowledgemets\rm
\end{center}
This work has been partially supported by Spanish Grant MTM2016-78751-P. The authors are most grateful for the constructive, detailed and useful remarks from an Associate Editor and an anonymous reviewer.

\end{document}